\documentclass{rspublic}
\usepackage{amsmath}
\usepackage{amsfonts}
\usepackage{amssymb}
\usepackage[Symbol]{upgreek}
\usepackage[nointegrals]{wasysym}
\usepackage{mathrsfs}
\usepackage{graphicx}
\usepackage{accents}
\newcommand{\real}{{\mathbb R}}
\newcommand{\bfa}{{\bf a}}
\newcommand{\bfb}{{\bf b}}
\newcommand{\bfu}{{\bf u}}
\newcommand{\bfv}{{\bf v}}
\newcommand{\bfw}{{\bf w}}
\newcommand{\bfx}{{\bf x}}
\newcommand{\bfy}{{\bf y}}
\newcommand{\bfz}{{\bf z}}
\newcommand{\bfg}{{\bf g}}
\newcommand{\bfV}{{\bf V}}
\newcommand{\bfU}{{\bf U}}
\newcommand{\bff}{{\bf f}}
\newcommand{\bfh}{{\bf h}}
\newcommand{\bfk}{{\bf k}}

\newcommand{\bfs}{{\bf s}}
\newcommand{\bbu}{\bar{\bf u}}
\newcommand{\bu}{\bar{u}}
\newcommand{\bp}{\bar{p}}
\newcommand{\baq}{\bar{q}}

\newcommand{\hbu}{\hat{\bf u}}
\newcommand{\hu}{\hat{u}}
\newcommand{\hv}{\hat{v}}
\newcommand{\hw}{\hat{w}}
\newcommand{\hp}{\hat{p}}
\newcommand{\hq}{\hat{q}}

\newcommand{\tp}{\tilde{p}}
\newcommand{\tbu}{\tilde{\bf u}}
\newcommand{\tu}{\tilde{u}}
\newcommand{\rr}{{{\mathbb R}^3}}
\newcommand{\Upo}{{\Omega}}
\newcommand{\bdy}{{\partial \Omega}}

\begin{document}

\title[Navier-Stokes Equations]{\large The Navier-Stokes equations in primitive variables}

\author[F. Lam]{F. Lam}

\affiliation{ }

\label{firstpage}

\maketitle

\begin{abstract}{Navier-Stokes Equations; Leray-Hopf Local Solution; Diffusion; Linearisation; Perturbation; Couette; Hagen-Poiseuille; Primitive Formulation}

The Navier-Stokes equations in the primitive formulation for incompressible flow describe the evolution of velocity and pressure, without recourse to vorticity. We show that, beyond the finite Leray-Hopf regularity interval, every postulated strong solution is accompanied by infinitely many diffusion-dominated percolations of arbitrary size, while the momentum deficit caused by the non-linearity is compensated by the pressure gradient. In the upper half space, we demonstrate how sequences of these collective companions can be re-scaled into an absurd singularity. Owning to the passive nature of the pressure, there exist no essential {\it a priori} bounds for establishing the uniqueness of primitive solutions. With the illustration of well-exploited examples of closed-form basic flows, we elucidate the reason why perturbations, infinitesimal or finite, instigate indeterminate states that render the concept of flow instability inadmissible. An effort has also been made to reappraise a number of important issues in fluid dynamics. Unfortunately, the primitive theory cannot serve as a reliable tool for prediction. 

Nevertheless, a dedicated effort has been made to elaborate {\it a priori} bounds for vorticity dynamics for ideal as well as real fluids. As a result, we are able to establish long-time regularity of the Cauchy problem for incompressible flows. A main conclusion is that no events of finite-time blow-up can ever occur in the Euler or Navier-Stokes equations for initial data of finite energy.
\end{abstract}

\tableofcontents
\section{Introduction}

The Navier-Stokes equations of motion are derived from the conservation principle of linear momentum, and the continuity (Navier 1823; Stokes 1845). For viscous incompressible flows, they are 
\begin{equation} \label{ns}
	\partial_t \bfu + (\bfu \cdot \nabla) \bfu  = - {\rho}^{-1} \nabla p + \nu \Delta \bfu,\;\;\; \nabla{\cdot}\bfu=0, 
\end{equation}
where the velocity vector $\bfu(\bfx,t)=(u,v,w)(\bfx,t)$, the coordinates $\bfx=(x_1,x_2,x_3)$, and the scalar $p(\bfx,t)$ denotes the pressure, and the kinematic viscosity, $\nu=\mu/\rho$, is the ratio of the dynamic viscosity and the density. We assume that body force has a potential which has been absorbed into the pressure. The dynamics is invariant in any homogeneous direction. In $\real^3$, we must observe
\begin{equation} \label{inv}
\bfu(\bfx,t)=\bfu(\bfx-\bfa,t),\;\;\;\nabla p(\bfx,t)=\nabla p(\bfx-\bfa,t),
\end{equation}
for arbitrary finite vector, $\bfa=(a_1,a_2,a_3)$. In cylindrical polar co-ordinates $(r,\theta,z)$, the axial $z$-direction is homogeneous in a number of specific cases. 

The initial solenoidal velocity is assumed to be smooth with compact support
\begin{equation} \label{ic}
 \bfu(\bfx,0)=\bfu_0(\bfx)\; \in \; C_c^{\infty}(\real^3). 
\end{equation} 
Specifically, initial data, $\bfu_0=(u_0,v_0,w_0)$, can be an element of the Schwartz class of functions (denoted by ${\mathcal S}$). Every component of $\bfu_0$ is a subset of $L^1(\real^3)$, defined to be smooth functions, where the derivatives of all orders decay faster than any polynomial,
\begin{equation} \label{swf}
\begin{aligned}
{\mathcal S}(\rr) = \big\{ v \in C^{\infty}(\rr): & \;\;\; \mbox{for all multi-indices}  \;\; \alpha, \;\beta \in {\mathbb N}, \;\; \mbox{finite}\;\; \bfx_0, \\
& \sup_{\bfx \in \rr}\big|(\bfx-\bfx_0)^{\beta} \; D^{\alpha} v(\bfx-\bfx_0)\big| < \infty \;\;\; \mbox{on} \;\; \rr \;\big\},
\end{aligned}
\end{equation}
where constant $\bfx_0$ is specified with the data. Any decay property of $v$ at large distances is with respect to $\bfx_0$, in view of the invariance (\ref{inv}).
In many applications, the requirement that initial data are elements of the Schwartz class ($\cal S$) can be unnecessarily restrictive. Instead, smooth functions of
\begin{equation} \label{ic2}
|\bfu_0| \leq \frac{B}{\big(A+|\bfx-\bfx_0|\big)^n}, \;\;\;(\mbox{finite}\; A > 0, \;n \geq 1),
\end{equation}
are admissible alternatives. The nominator $B$ is a constant or specific bounded functions, $|B(\bfx,\bfx_0)| \leq 1$, for instance. Initial vorticity is found from 
\begin{equation} \label{vtic}
\upomega_0(\bfx)=\nabla{\times}\bfu_0(\bfx).
\end{equation}

For any arbitrary smooth function, $\bff=\bff(t)$, $\bff(0)=\bff'(0)=0$, spontaneous solutions to the Navier-Stokes or the Euler equations are given by $(\bfu,\nabla p)=(\bff,- \rho \bff')$. In order to avoid under-determinacy of the Cauchy problem, we need to impose a `boundary condition' at infinity or an equivalent constraint. Physics suggests that the energy must be finite,
\begin{equation} \label{c2}
\int_{\rr} |\bfu|^2\; \rd \bfx < \infty,
\end{equation}
for $t \in \Upo_T$, where $\Upo_T=[0, T< \infty)$. In the Sobolev space, the finiteness translates into $u_i \in W_0^{0,2}(\rr)$ or $\bfu$ having compact support. Instead of (\ref{c2}), we require that the velocity tends to zero at large distances,
\begin{equation} \label{vtbc}
\bfu \rightarrow 0 \;\;\;\mbox{as} \;\;\; |\bfx| \rightarrow \infty.
\end{equation}
At this stage, we can only state the decay loosely, as we do not have more knowledge about the precise decay rate. In subsequent discussions, all similar expressions of attenuation are understood in the relative sense, $|\bfx-\bfx_0| \rightarrow \infty$. Conversely, this property implicitly expresses a bound that $\bfu \in W_0^{0,1}$ on $\rr$. In an extreme scenario, for a source singularity in potential flow, it is known that the velocity decays at infinity $\sim O(|\bfx|^{-2})$. 
By virtue of Gauss's divergence theorem, we have 
\begin{equation} \label{dv}
\int_{\rr} \nabla{\cdot}\bfu\; \rd \bfx = \lim_{r\rightarrow \infty}\int_S \bfu \cdot \vec{n} \; \rd \bfs = 0,
\end{equation}
where, $S=4 \pi r^2$, stands for the surface area of a ball centred at the origin $\bfx_0$, the radius $r$ of which tends to large values. It follows that, as one possibility, the continuity is valid provided $\bfu \sim O(|\bfx|^{-2})$ as $|\bfx|\rightarrow \infty$.\footnote{Condition (\ref{dv}) merely states that the sum of the radial velocity must be zero. There remains a possible scenario that some fluid enters one portion of the spherical surface while the exact same amount leaves at another portion. Moreover, we do not have complete knowledge about local tangential components; the whole ball might be under a uniform rotation, for instance. As $r$ tends to be unbounded, these non-zero velocities cannot be sustained for finite energy.} Similarly, the solenoidal vorticity implies that, at least, $\upomega \sim O(|\bfx|^{-2})$ at infinity. Hence `boundary condition' (\ref{vtbc}) does not compromise analytical scopes of flow-fields consisting of $\bfu$-$\upomega$ and, in fact, it serves our purposes well, compared with periodic flows.

In bounded domains with $C^2$ boundaries, we relax the data 
\begin{equation} 
\bfu(\bfx,0)=\bfu_0(\bfx)\; \in \; C^{\infty}(\Upo),
\end{equation}
as we are dealing with flows of finite energy. To be definitive, we also impose the no-slip condition, $\bfu(\bfx,t)=0,\forall \bfx \in \bdy$.
Taking divergence on (\ref{ns}), and using the continuity, we obtain the Poisson equation for the pressure 
\begin{equation} \label{ppoi}
	\Delta p(\bfx)= - \rho \: \big(\: |\nabla \bfu|^2 - |\nabla{\times}\bfu|^2  \:\big)(\bfx).
\end{equation}
In bounded domains, this is a Neumann boundary value problem.

The conservation laws are universal and, indeed, they do not depend on the measuring units in any particular problem. We view the equations of motion (\ref{ns}) (and associated initial and boundary conditions) as inherently dimension-independent with respect to the SI standard references of length, time, and mass. To be consistent, the viscosity is a physical property of fluids, hence it is made non-dimensional in reference to unit kinematic viscosity $\nu_0{=}1 \:([\:\texttt{m}^2 \:\texttt{s}^{-1}\:])$. In fact, a reference kinematic viscosity rather than unity can be chosen, as long as the ambient parameters are standardised, for instance, at a fixed set of the SI temperature and density.

The theoretical framework, (\ref{ns}) to (\ref{ppoi}), is known as the initial-boundary value problem in the primitive variables, ($\bfu,p$). In a nutshell, the primitive equations admit local-in-time smooth solutions, depending on the size of  the initial data (Leray 1934$a$; Hopf 1951); the global regularity is out of the question, due to lack of suitable {\it a prioir} bounds. The popularly-adopted classification of two distinct modes of flow, streamlined laminar and fluctuating turbulent, is completely based on the primitive setting (Reynolds 1883). Perturbation approaches on shear flow instabilities, linear or non-linear, are all formalised in terms of the $\bfu{-}p$ variables. Assisted by modern computers and algorithms, numerical solutions are readily available in applications, despite the absence of rigorous mathematical justifications. It has been well-known that the computations themselves often stimulate new analytical problems, such as error-controls in discretisation, convergence and stability of numerics. In the present note, we examine the compatibility and the global regularity of the primitive formulation. 

By examining Helmholtz's decomposition of vector fields on $\rr$, we establish the vorticity-velocity compatibility,  $\Delta \bfu = - \nabla{\times}\upomega$. The importance lies in the fact that fluid dynamics is tantamount to solving vorticity equation (see equation (\ref{vort})) in conjunction with this kinematic constraint. Then, we make an effort on derivation of {\it a priori} bounds for the vorticity equation. These estimates are crucial in determining the global regularity of the Cauchy problem.

\section{Non-uniqueness}

Let the Leray-Hopf regular solution be $(\bfu^*,p^*)$, and its existence time interval $0 < t \leq t^*$. 
Suppose that there is an extended local solution, or even a global solution of the Navier-Stokes equations, denoted by $(\bfu,p)(\bfx, t > t^*)$ with the starting data (\ref{ic}). We consider the case that the solution is smooth, or sufficiently differentiable. The extended solution may be obtained by either an analytic construction or a numerical procedure. We are particularly interested in the numerical solutions established from a projection method, or a fractional step technique, where the primitive variables, $\bfu(\bfx,t)$ and $p(\bfx,t)$, are iterated in parallel. It is evident that the solutions $(\bfu,p)$ and $(\bfu^*,p^*)$ must coincide over the interval $0 \leq t \leq t^*$. 

\subsection{Weak formulation and buffer functions}

All the derivatives in the equations of motion are interpreted in the sense of distributions. Let the test functions,
\begin{equation*}
\upvarphi(\bfx,t)=(\varphi_1,\;\varphi_2,\;\varphi_3)(\bfx,t),\;\; (\in C^{\infty}_c(\real^3)\times\real)
\end{equation*}
be vectors. Multiplying the first equation of (\ref{ns}) by $\upvarphi$, integrating the result over space and time, followed by integration by parts, 
we obtain the weak formulation
\begin{equation} \label{wkns}
\int_{\real}\!\int_{\real^3} \bfu \: \Big(\: \frac{\partial \upvarphi}{\partial t} + \nu \Delta \upvarphi\:\Big) \rd \bfx \rd t = \sum_{i,j=1}^3\int_{\real}\!\int_{\real^3} u_i \frac{\partial u_j}{\partial x_i} \varphi_j \: \rd \bfx \rd t - \int_{\real}\!\int_{\real^3} p\:\nabla.\upvarphi \: \rd \bfx \rd t.
\end{equation}
This is Leray's definition of weak solutions (Leray 1934$a$), where the velocity $\bfu$ is square-integrable, and the pressure $p$ locally integrable. If the flow field is found to be incompressible and sufficiently regular, the continuity (\ref{ns}) makes sense. Leray also showed that weak solutions (\ref{wkns}) satisfy energy inequality,
\begin{equation} \label{ienergy} 
\frac{1}{2} \int_{\real^3} |\bfu(\bfx,t)|^2 \rd \bfx + \nu \int_0^t\!\!\int_{\real^3} |\nabla \bfu(\bfx,t)|^2 \rd \bfx \rd t \; \leq \; \frac{1}{2} \int_{\real^3} |\bfu(\bfx,0)|^2 \rd \bfx.
\end{equation}

As proposed by Hopf (1951), we may emphasise the role of the initial condition,  as long as the test functions are re-defined as
\begin{equation*}
\upvarphi(\bfx,t)=(\varphi_1,\;\varphi_2,\;\varphi_3)(\bfx,t),\;\; (\in C^{\infty}_c(\real^3)\times [0,\infty))
\end{equation*}
on a set of measure zero (time-wise). The formulation (\ref{wkns}) is extended to 
\begin{equation} \label{wknsp}
\begin{split}
\int_0^{\infty}\!\!\!\int_{\real^3} \bfu \: \Big(\: \frac{\partial \upvarphi}{\partial t} &+ \nu \Delta \upvarphi\:\Big) \rd \bfx \rd t + \int_{\real^3 }\bfu_0 \:\upvarphi(0) \: \rd \bfx   \\
& - \sum_{i,j=1}^3\int_0^{\infty}\!\!\!\int_{\real^3} u_i \frac{\partial u_j}{\partial x_i} \varphi_j \: \rd \bfx \rd t + \int_0^{\infty}\!\!\!\int_{\real^3} p\:\nabla{\cdot}\upvarphi \: \rd \bfx \rd t=0.
\end{split}
\end{equation}
Apparently, the revised formula does not contain noticeable differences compared to Leray's. However, given every initial velocity $\nabla.\bfu_0=0$, the linear kernel, $\partial_t - \nu \Delta $, is data-preserving in the sense that the solenoidal field is maintained throughout the motion. Thus the continuity of (\ref{ns}) makes sense and is readily satisfied. This observation suggests that the formulation (\ref{wknsp}) can be interpreted as a composition of two degenerated dynamics,
\begin{equation} \label{wkns1}
\int_0^{\infty}\!\!\!\int_{\real^3} \bfu \: \Big(\: \frac{\partial \upvarphi}{\partial t} + \nu \Delta \upvarphi\:\Big) \rd \bfx \rd t + \int_{\real^3 }\bfu_0 \:\upvarphi(0) \: \rd \bfx  = \int_0^{\infty}\!\!\!\int_{\real^3} \bfb \upvarphi \: \rd \bfx \rd t,
\end{equation}
and
\begin{equation} \label{wkns2}
\sum_{i,j=1}^3\int_0^{\infty}\!\!\!\int_{\real^3} u_i \frac{\partial u_j}{\partial x_i} \varphi_j \: \rd \bfx \rd t - \int_0^{\infty}\!\!\!\int_{\real^3} p\:\nabla{\cdot}\upvarphi \: \rd \bfx \rd t=\int_0^{\infty}\!\!\!\int_{\real^3} \bfb \upvarphi \: \rd \bfx \rd t,
\end{equation}
if there exists a function $\bfb=\bfb(\bfx,t)$ mediating the parities. Although the Cauchy problem (\ref{wkns}) or (\ref{wknsp}) shares identical initial data with formulation (\ref{wkns1}), their solutions do not overlap in general, except over a very short interval $[0,\varepsilon]$, as the former setting aims at wholly non-local and non-linear solutions. 

To determine such a function $\bfb$, we first state the precise definition of the solutions of (\ref{ns}). It is known that the regular Leray-Hopf solutions do exist over a finite period of time. As strong solutions are also weak solutions, we shall work with regular functions and tactically put aside the regularity issue which entails {\it a priori} bounds.  We call ($\bfu,p$) a Navier-Stokes solution if the pair satisfy all of the following conditions: 
\begin{equation} \label{nsol}
\left.
\begin{aligned}
& 1.\;\; \mbox{The velocity field is solenoidal}, \nabla{\cdot}\bfu=0; \;\mbox{at}\; t > 0, \\
& \hspace{5mm} \mbox{it depends continuously on its initial solenoidal data} \;\bfu_0; \\
& 2.\;\; \mbox{Linear momenta are preserved at every instant} \;t > 0;\\
& 3.\;\; \mbox{The pressure is determined by Poisson's equation},\\
& \hspace{10mm}\Delta p = - \rho \nabla{\cdot}\big(\: (\bfu.\nabla)\bfu \:\big) \;\;
(\mbox{up to arbitrary function} \;p_0(t));\\
& 4.\;\; \mbox{The law of energy conservation holds}, \;\mbox{for}\;\nu>0,\\
& \hspace{10mm}\frac{1}{2} \int_{\real^3} |\bfu|^2 \rd \bfx + \nu \int_0^t\!\!\int_{\real^3} |\nabla \bfu|^2 \rd \bfx \rd t \; = \; \frac{1}{2} \int_{\real^3} |\bfu_0|^2 \rd \bfx.
\end{aligned}
\hspace{10mm} \right\}
\end{equation}
In addition, the Neumann boundary value problem for the pressure is self-consistent. Conceptually, the last two constraints are {\it de facto} the derivatives of the first two.

Let $\bfb \in C^{\infty}_c(\real^3)$ or $\bfb \in C^{\infty}(\Upo)$ in bounded domains. The second condition implies that the momentum equation in (\ref{ns}) or its weak forms, (\ref{wkns}) and (\ref{wknsp})-(\ref{wkns2}), can be interpreted as one of the following equivalent statements:
\begin{equation} \label{efs}
\begin{aligned}
\partial_t \bfu &+ (\bfu \cdot \nabla) \bfu  = - {\rho}^{-1} \nabla p + \nu \Delta \bfu \;\;\;\; &\Longleftrightarrow \\
\partial_t \bfu &- \nu \Delta \bfu= - (\bfu \cdot \nabla) \bfu - {\rho}^{-1} \nabla p = \bfb \;\;\;\; &\Longleftrightarrow\\
\partial_t \bfu &- \nu \Delta \bfu - \bfb = - (\bfu \cdot \nabla) \bfu - {\rho}^{-1} \nabla p - \bfb \;\;\;\; &\Longleftrightarrow\\
\partial_t \bfu &+ (\bfu \cdot \nabla) \bfu  + {\rho}^{-1} \nabla p - \nu \Delta \bfu=0,
\end{aligned}
\end{equation}
where each form respects the conservation principle for every solenoidal field $\bfu$. We single out the `forced' diffusion dynamics, 
\begin{equation*} 
\partial_t \bfu - \nu \Delta \bfu = \bfb,
\end{equation*}
subject to the solenoidal initial data $\bfu_0$. Its solution can be expressed as
\begin{equation*} 
\begin{split}
\bfu(\bfx,t) & = \int_{\real^3} H(\bfx{-}\bfy,t)\bfu_0(\bfy) \rd \bfy + \int_0^t \!\int_{\real^3} H(\bfx{-}\bfy,t{-}t')\bfb(\bfy,t') \: \rd \bfy \rd t' \\
\quad & \hspace{5mm} = \int_{\real^3} H(\bfy,t)\bfu_0(\bfx{+}\bfy) \rd \bfx + \int_0^t \!\int_{\real^3} H(\bfy,t{-}t')\bfb(\bfx{+}\bfy,t') \: \rd \bfy \rd t',
\end{split}
\end{equation*}
where $H$ is the free-space heat kernel,
\begin{equation} \label{heat}
H(\bfx{-}\bfy,t)=\big(\:4 \pi \nu t\:\big)^{-3/2} \exp \big( - {|\; \bfx-\bfy\; |^2}/(4 \nu t)\:\big).
\end{equation}
It is this property of continuous dependence on initial data that motivates us to split the pressure-related momentum equation into two separate parts. From (\ref{efs}), we see that the pressure can be found from the Neumann problem
\begin{equation} \label{ppoi0}
\Delta p = - \rho \nabla{\cdot}\big(\: (\bfu{\cdot}\nabla)\bfu + \bfb \:\big),
\end{equation}
for a proper choice of the function $\bfb$, because the non-linear term is known, though the convection now becomes moderated. To meet every criterion of (\ref{nsol}), we specify that the buffer function is calorically mollifiable, and  within free multiples,
\begin{equation} \label{buff}
\bfb(\bfx,t)=\upgamma(t) \uppsi(\bfx),
\end{equation}
where time-dependent $\upgamma$ is any regular function, and $\uppsi$ satisfies one of the following conditions:
\begin{enumerate}
\item $\nabla{\cdot}\uppsi=0$, and $\bfu{\cdot}\uppsi=0$ (orthogonal for energy conservation);
\item $\uppsi=\nabla{\times}\upphi$, and $\bfu{\cdot}\uppsi=0$;
\item $\uppsi=\bfv_1{\times}\bfv_2$, $\bfv_1$ and $\bfv_2$ being irrotational, and $\bfu{\cdot}\uppsi=0$;
\item $\uppsi=(\nabla\psi) {\times}(\nabla \phi)$, for scalars $\psi$ and $\phi$, and $\bfu{\cdot}\uppsi=0$;
\item $\uppsi=\Delta \upchi$, for solenoidal $\upchi$, and $\nabla \bfu{\cdot} \nabla \upchi=0$;
\item $\uppsi = \nabla \phi(\bfx)$, and $\phi$ is harmonic.
\end{enumerate}
Choosing a buffer from the list, we see that the Poisson's equation (\ref{ppoi0}) is reduced to the third condition of (\ref{nsol}). To balance the momentum, we simply set
\begin{equation} \label{dpp0}
\nabla p = - \rho \big(\: (\bfu{\cdot}\nabla)\bfu + \bfb \:\big),
\end{equation}
where the right-hand side is completely given. For bounded domains, the evaluation is slightly more complicated but is tractable; the dependence of the pressure on $\bfb$ can be found from the Neumann boundary conditions, unless we further impose the no-slip on $\bfb$. These relations rectify the auxiliary nature of the pressure, or more precisely its gradient $\nabla p$. 

We point out that list (\ref{buff}) is by no means comprehensive. Linear superpositions of these classes are also good candidates. Specifically, the Helmholtz's vector decomposition is included. The last category represents weak solutions of (\ref{wknsp}) for any solenoidal test function (Serrin 1962). 

\subsection{Percolation and pressure compensation}

Now, we return to the issues raised at the beginning of the present section. Given a regular solution $(\bfu,p)$ of (\ref{wkns}), or (\ref{ns}), which has been postulated to exist beyond the Leray-Hopf time $t>t^*$. Consider the following superimposed pair
\begin{equation} \label{ss}
\big(\: \bfu+\alpha \hbu,\;\;p+ \hp/\beta \: \big)(\bfx,t),
\end{equation}
where the parameters, $\alpha$ and $\beta=\beta(t)$, are arbitrary but finite. The pair ($\hbu,\hp$) are regarded as perturbations to the postulated solution. We substitute the superposition (\ref{ss}) into the Navier-Stokes equations (\ref{ns}). Since solution $(\bfu,p)$ satisfies the equations of motion, we are left with
\begin{equation} \label{nss}
\partial_t \hbu + \alpha (\hbu{\cdot}\nabla)\hbu + (\bfu{\cdot}\nabla)\hbu + (\hbu{\cdot}\nabla)\bfu = -(\alpha \beta \rho)^{-1} \nabla \hp + \nu \Delta \hbu,\;\;\;\nabla{\cdot}\hbu=0.
\end{equation}
To solve these exactly, we characterise the following diffusion problem as part of the solution:
\begin{equation} \label{dfu}
\partial_t \hbu - \nu \Delta \hbu = \bfb,
\end{equation}
subject to the initial conditions
\begin{equation*}
\hbu_0(\bfx)=\bfu^*(\bfx,t^*)=\bfu(\bfx,t^*).
\end{equation*}
An extension of the Leray-Hopf solution is feasible, as the diffusion equation has solution 
\begin{equation} \label{ds}
\hbu(\bfx,t)=\int_{\real^3} H(\bfy,t) \: \hbu_0(\bfx{+}\bfy) \: \rd \bfy + \int_{t^*}^t \!\int_{\real^3} H(\bfy,t{-}t')\bfb(\bfx{+}\bfy,t') \: \rd \bfy \rd t'.
\end{equation}
Since the initial data $\nabla{\cdot}\bfu(\bfx,t^*)=\nabla{\cdot}\bfu^*=0$, the solution $\hbu$ remains solenoidal for all times, as indicated by the expression in (\ref{ds}), and the properties of the buffers. Moreover, $\hbu$ is globally regular in space and in time, by virtue of the well-known smoothing properties of the heat kernel. In fact, we only need 
\begin{equation*}
\bfu_0(\bfx,t^*) \in L^1_{loc}(\bfx)\; \forall \bfx \in \real^3.
\end{equation*}

The analytical structure of equation (\ref{nss}) suggests that we can compensate the momentum deficiency due to $\hbu(\bfx,t)$ exactly, by a complement of pressure gradient $\nabla\hp$. Analytically, equations (\ref{nss}) are fully satisfied as long as the pressure $\hp$ is found from
\begin{equation} \label{bns}
\nabla \hp= - \; A_0 \; \big( \: \alpha (\hbu{\cdot}\nabla)\hbu +(\bfu{\cdot}\nabla)\hbu + (\hbu{\cdot}\nabla)\bfu \: \big) - \rho \beta \bfb = - \bfg,
\end{equation}
where $A_0=\alpha \beta \rho $. In practice, only the gradient rather than the actual pressure matters. Nevertheless, we proceed to describe a general approach. Let us introduce the abbreviations:
\begin{equation*}
\bfV[\bfu_x]=(\partial u/\partial x \;\;\;\partial v/\partial x\;\;\;\partial w/\partial x)',\;\;\;\;\;\;
S[u]=(\partial u/\partial x\;\;\; \partial u/\partial y\;\;\;\partial u/\partial z),
\end{equation*}
where the symbol $\bfV$ stands for a $3{\times}1$ column matrix operating on vector $\bfu$, $S$ for a $1{\times}3$ row matrix on a scalar. Taking divergence on equation (\ref{bns}), or equivalently on (\ref{nss}), we readily establish the following equation:
\begin{equation} \label{bsp}
\begin{split}
	\Delta \hp(\bfx)= {-} A_0\: \Big(\: & \alpha \big(\: S[\hat{u}]\:\bfV[\hbu_x]+S[\hv]\:\bfV[\hbu_y]+S[\hw]\:\bfV[\hbu_z] \: \big) \; + \\
	\quad & 2 \big( \: S[u]\: \bfV[\hbu_x]+S[v]\: \bfV[\hbu_y]+S[w]\: \bfV[\hbu_z] \: \big)\:\Big)(\bfx) = -F(\bfx).
\end{split}
\end{equation}
Evidently, the pressure $\hp$ is determined, up to a constant, for fixed $\alpha$ and $\beta$ because each of the $18$ velocity terms is known, or
\begin{equation} \label{psol}
	\hp = \frac{1}{4 \pi}\int_{\real^3} \frac{1}{|\bfx-\bfy|} F(\bfy) \: \rd \bfy + \hp_0.
\end{equation}
For every postulation $(\bfu,p)$, the full Navier-Stokes system (\ref{ns}) admits the companion solutions (\ref{ss}) over the time interval $t > t^*$. We call the pair $(\hbu,\hp)$ a (viscous) percolation, because they diffuse everywhere at all times.  

As a well-known practice, the energy conservation is obtained by integrating complete momentum (\ref{nss}) or the diffusion (\ref{dfu})  
\begin{equation} \label{energy}
\frac{1}{2} \int_{\real^3} |\hbu|^2 \rd \bfx + \nu \int_0^t\!\!\int_{\real^3} |\nabla \hbu|^2 \rd \bfx \rd t= \frac{1}{2} \int_{\real^3} |\hbu_0|^2 \rd \bfx,
\end{equation}
which is scale-invariant to parameter $\alpha$. We draw our attention to the converse: if we rely on the law as a bound, we cannot distinguish a percolation from a complete non-linear solution. In particular, the percolation is arbitrary, while the full solution, or a result of numerical computation, is expected to be unique.

In numerical schemes, such as a fractional step algorithm or a projection method, the appearance of non-zero buffer function $\bfb$ is inevitable. Over the mesh and time, discretisation errors of the operator $\partial_t {-} \nu \Delta$ may well be in the same order of those on the non-linearity and the pressure (cf. (\ref{efs})), implying the existence of buffers. If we strictly solve the primitive equations, we have no reliable means to detect or filter out the differences of the errors which evolve with the computations. In the case of zero buffer, the structure of the primitive formulation still allows the flow-field partition, according to the criteria (\ref{nsol}).

It has been suggested that, the pressure plays the role of a Lagrangian multiplier which ensures incompressibility. This is a misleading view, given the partition. The gradient is determined by (\ref{bns}), in a passive manner, resulting in the superficial Poisson equation (\ref{ppoi}). One implication is that use of the Helmholtz-Leray projection, which eliminates the pressure in conjunction with the continuity, lacks generality.

\subsection{Multiplicity inside a cube}

Let the domain be the cubic box with all sides of unity, 
$\Upo_b: \; 0 \leq x,y,z \leq 1$.
We impose the no-slip condition on the walls of the cube.\footnote{It has been known that the initial-boundary value problems for incompressible flows ($\nu \geq 0$) are ill-posed in periodic domains with periodic boundary conditions (Lam 2018$b$). We emphasise the fact that the conjecture of anomalous energy dissipation for inviscid-limit continuum turbulence by Onsager (1949) was made on an assumption of the periodicity for the primitive variables. In the last paragraph on p.286, Onsager discussed his hypothesised breakdown of energy conservation in terms of the Fourier formulation of the governing equations on the macroscopic scales (cf. his equations (15) and (17) of p.284).} We choose $\bfb=0$ so that the velocity is defined by pure diffusion
\begin{equation} \label{dfu0}
\partial_t \hbu - \nu \Delta \hbu = 0.
\end{equation}
By analogy, the companion velocity to any extended Leray-Hopf flow is still evaluated by the first integral of (\ref{ds}). Instead of (\ref{heat}), the heat kernel now reads
\begin{equation} \label{cu}
	H_f(x,x',t) \times H_f(y,y',t) \times H_f(z,z',t),
\end{equation}
where the finite-line heat kernel has the form
\begin{equation} \label{hf} 
\begin{split}
	H_f(x,x',t) & =2 \;\sum_{k=1}^{\infty}\;\sin(k\pi x)\sin(k\pi x')\exp\big(-\nu k^2 \pi^2  t \big)\\
	\quad & = \frac{1}{\sqrt{4 \pi \nu t}}\sum_{k=-\infty}^{\infty}\Big[\: \exp\Big( {-}\frac{(x{-}x'{-}2k)}{4 \nu t}\Big){-}\exp\Big({-}\frac{(x{-}x'{+}2k)}{4 \nu t} \Big)\:\Big].
	\end{split}
\end{equation}
Since $\hbu = \bfu = 0$ on the walls, the zero Neumann boundary condition is found from (\ref{bns}), namely $\partial \hp/\partial \vec{n}|_{\bdy_b}=0$ in the present case. The pressure is governed by a compatible system because, from (\ref{bsp}) and (\ref{bns}), we have
\begin{equation*}
	\int_{\Upo_b} F(\bfx) \rd \bfx = - \int_{\Upo_b} \Delta \hp (\bfx) \: \rd \bfx = \int_{\Upo_b} (\nabla{\cdot} \bfg)(\bfx) \: \rd \bfx = \int_{\bdy_b} (\bfg{\cdot}\vec{n}) \: \rd s=0,
\end{equation*}
by virtue of Gauss' divergence theorem. For the demonstration of the existence, it is preferable to express the pressure solution (\ref{psol}) in terms of the Neumann function, $N_b$, which is defined by 
\begin{equation} \label{dnc}
\Delta N_b=-\delta(\bfx,\bfx')+1,\;\;\;\partial N_b/\partial \vec{n}\big|_{\bdy_b}=0,
\end{equation}
where $\delta=\delta(x)$ is Dirac's delta. We start from a cosine expansion for the delta in one dimension,
\begin{equation*}
\delta(x-x')=\frac{1}{L}+\frac{2}{L}\sum_{l=1}^{\infty} \cos\Big(\frac{l\pi x}{L}\Big)\cos\Big(\frac{l\pi x'}{L}\Big) = \frac{1}{L}+\frac{2}{L}\sum_{l} c_l^x c_l^{x'},
\end{equation*}
for $0 < x,x'< L$. Then it is easy to establish, with justifications, see, for instance, Chapter 9 of Roach (1982), 
\begin{equation*}
\begin{split}
-\delta(\bfx-\bfx')+\frac{1}{L^3}&=-\frac{2}{L^3}\bigg(\sum_{l} c_l^x c_l^{x'}+\sum_{m} c_m^y c_m^{y'}+\sum_{n} c_n^z c_n^{z'}+2\sum_{l,m} c_l^x c_l^{x'}c_m^y c_m^{y'}\\
\quad & + 2 \sum_{m,n} c_m^y c_m^{y'}c_n^z c_n^{z'}
+2 \sum_{n,l} c_n^z c_n^{z'}c_l^x c_l^{x'} + 4\sum_{l,m,n} c_l^x c_l^{x'}c_m^y c_m^{y'}c_n^z c_n^{z'}\bigg).
\end{split}
\end{equation*}
By symmetry, the Neumann function ($L=1$) must have the neat expansion,
\begin{equation*}
N_b(\bfx,\bfx')=8\sum_{l,m,n=0} A(l,m,n,x',y',z')\; c_l^x c_m^y c_n^z,
\end{equation*}
which ensures the derivative homogeneous data on the boundary. From the differential equation (\ref{dnc}), we find the required Fourier expansions in cosines
\begin{equation*} 
\begin{split}
	N_c(\bfx,\bfx')&=-\frac{8}{\pi^2} \sum_{l,m,n=0}^{\infty}\!\Big(\:A_{lmn}\cos(l \pi x)\cos(l \pi x')\: \cos(m \pi y)\cos(m \pi y')\\
	\quad & \hspace{35mm} \cos(n \pi z)\cos(n \pi z')\:\Big) \Big/ \Big( l^2+m^2+n^2 \Big),
\end{split}
\end{equation*}
where the normalising factors are
\begin{equation*}
\begin{split}
A_{000}&=0; \;\; A_{00n}=A_{0m0}=A_{l00}=1/4; \\
A_{0mn}&=A_{l0n}=A_{lm0}=1/2; \;\;\mbox{and}\;\; A_{lmn}=1, \;\;\mbox{for}\;\; l,m,n\geq1.
\end{split}
\end{equation*}

An alternative way to determine the pressure is as follows. The Neumann boundary value problem,
\begin{equation*}
\Delta \hp = - F,\;\;\;\;\;\nabla \hp\big|_{\bdy_b}=0,
\end{equation*}
can be reduced to a Dirichlet problem, as $-\nabla \hp|_{\partial \Upo_b}$ is found from momentum (\ref{ns}). This approach may have an advantage in the circumstances where we have computed the velocities as well as their derivatives. By virtue of (\ref{bsp}), what we need to analyse is a vector Poisson equation for $\bfg$ with homogeneous boundary data
\begin{equation} \label{pdv}
-\Delta(\nabla\hp)=\nabla F,\;\;\;\;\;\bfg\big|_{\bdy_b}=-\nabla\hp\big|_{\bdy_b}=0.
\end{equation}
Solving (\ref{pdv}) without specifying the unknown gauge, we find the gradient 
\begin{equation*} 
\nabla \hp(\bfx) = \int_{\Upo_b} G_b(\bfx,\bfx')\nabla F (\bfx') \: \rd \bfx',
\end{equation*}
where the Green's function $G_b$ equals
\begin{equation*}
	\frac{8}{\pi^2} \!\!\sum_{l,m,n=1}^{\infty}\!\!\sin(l \pi x)\sin(l \pi x') \sin(m \pi y)\sin(m \pi y')\sin(n \pi z)\sin(n \pi z') /( l^2{+}m^2{+}n^2 ).
\end{equation*}

Our partition proposition, as afforded by (\ref{wkns1}), (\ref{wkns2}) or (\ref{efs}), tactically bypasses the task of handling the non-linearity. It is instructive to recall that Stokes' system,
\begin{equation*} 
\nu \Delta \bfu - \nabla p /\rho = \sigma \bfu;\;\;\;\nabla{\cdot}\bfu=0,\;\;\;(\bfu|_{\Upo_b}=0)
\end{equation*}
has been extensively analysed for the spectra of eigen-values ${\sigma_k}$, and the normalised eigen-functions. A weak solution of (\ref{wkns}) is represented, in a Faedo-Galerkin approximation, by the complete orthonormal basis, in analogous ways to Fourier series expansions. Technically, this is a linear strategy as well. However, it by no means inspires the `lower-order' view that prioritises the dynamic role of viscous diffusion.

\subsection{Artificial singularity in upper half space}

Denote the flow domain by $\Upo_h: -\infty < x,y < \infty;\;0 \leq z < \infty$. Let $(\bfu^*,p^*)$ be a Leray-Hopf regular solution in $0\leq t\leq t^*$, where $\bfu^*(x,y,0,t)=0$. Suppose that there is a mean solution, ($\bbu,\bp$), which extends $(\bfu^*,p^*)$ beyond $t^*$. At $t=t^*$, these two solutions coincide. The no-slip on the plane $z=0$ holds for the mean motion $t>t^*$. Consider a bounded time interval $[t^*\leq t \leq T]$, which $T$ lies within the existing time of the mean motion. We then sub-divide it into $n$ equal parts, $\Delta t= (T-t^*)/n$. Denote $t_k=t^*+k \Delta$ or
\begin{equation*}
t^*< t^* + \Delta t <  t^* + 2\Delta t < \cdots < t^* + k \Delta t < \cdots < T.
\end{equation*}

We introduce the first percolation ($\hbu^{(1)},\hp^{(1)}$) on the mean motion over $[t^*,t_1]$, starting from
\begin{equation*}
\hbu_0(\bfx) = \bfu^*(\bfx,t^*)=\bbu(\bfx,t^*).
\end{equation*}
The solution of (\ref{dfu0}) is expressed in the integral
\begin{equation} \label{hsp}
\hbu^{(1)}(\bfx,t)=\int_{\Upo_h} H_s(\bfx,\bfx',t) \hbu_0(\bfx') \: \rd \bfx',\;\;\;(t^* < t < t_1)
\end{equation}
where the kernel $H_s$ is familiar and vanishes on the plane $z=0$,
\begin{equation} \label{hhsp}
H_s(\bfx,\bfx',t)=H_1(x,x',t) \times H_1(y,y',t)\times \big(H_1(z,z',t)-H_1(z,-z',t)\big),
\end{equation}
and the function $H_1$ is the free-line heat kernel
\begin{equation} \label{h1}
	H_1(x,x',t)={1}/{\sqrt{4 \pi \nu t}}\; \exp\big(\:- {(x-x')^2}/{(4 \nu t)}\:\big).
\end{equation}
Given the no-slip on $z=0$, the Neumann function can be found by the method of images. The formula for the solution of (\ref{bsp}) reads
\begin{equation} \label{hspp}
\hp^{(1)}(\bfx)= \frac{1}{4 \pi}\int_{\Upo_h} \Big(\:\frac{1}{r_-}+\frac{1}{r_+}\:\Big)F^{(1)} (\bfx') \: \rd \bfx' + \hp_0,
\end{equation}
where $F=F^{(1)}$ contains the first percolation $\hbu^{(1)}$, and
\begin{equation*}
r_{\mp}(\bfx,\bfx')=\Big(\:(x-x')^2+(y-y')^2+(z{\mp}z')^2\:\Big)^{1/2}.
\end{equation*}
The first companion solution is
\begin{equation} \label{cs1}
(\bfu_c, p_c) = \big(\bbu+\alpha^{(1)}\hbu^{(1)},\;\bp+\hp^{(1)}/\beta^{(1)}\big)
\end{equation}
over $[t^*,t_1]$. Next we introduce a second percolation, ($\hbu^{(2)},\hp^{(2)}$) on $(\bfu_c,p_c)$ over $[t_1,t_2]$ with starting data
\begin{equation*}
\hbu_0(\bfx)=\bfu_c(\bfx,t_1).
\end{equation*}
We repeat our analysis to fix the second percolation ($\hbu^{(2)},\hp^{(2)}$), and renew the forcing function $F^{(2)}=F^{(2)}(\hbu^{(2)})$ in (\ref{hspp}) to settle $\hp^{(2)}$. Now the companion is updated to
\begin{equation} \label{cs2}
(\bfu_c, p_c) = \big(\bbu+\alpha^{(1)}\hbu^{(1)}+\alpha^{(2)}\hbu^{(2)},\;\bp+\hp^{(1)}/\beta^{(1)}+\hp^{(2)}/\beta^{(2)}\big)
\end{equation}
over $[t^*,t_2]$. Assume that, after $k$-time percolations, we obtain the sum
\begin{equation} \label{csk}
\bbu+\sum_{i=1}^k\alpha^{(i)}\hbu^{(i)},\;\bp+\sum_{i=1}^k\hp^{(i)}/\beta^{(i)} 
\end{equation}
over $[t^*,t_k]$. To determine the $(k+1)$th percolation, we substitute $(\bfu_c+\alpha^{(k+1)}\hbu^{(k+1)}$, $p_c+\hp^{(k+1)}/\beta^{(k+1)})$ into the Navier-Stokes and subtract the companion, we obtain the analogous equations (\ref{nss}) and (\ref{bsp}) for $\hbu^{(k+1)}$ and $\hp^{(k+1)}$ respectively. We solve 
\begin{equation*}
\partial_t\hbu^{(k+1)}-\Delta\hbu^{(k+1)}=0,
\end{equation*}
subject to data $\hbu_0=\bfu_c(\bfx,t_k)$, and its solution is expressed in the integral (\ref{hsp}). The pressure, $\hp^{(k+1)}$, is given in the formula (\ref{hspp}) with $F^{(k+1)}=F^{(k+1)}(\hbu^{(k+1)})$. Evidently, all of the the percolations are exact and smooth. Specifically, the companion solution to ($\bbu,\bp$) over $[t_{k-1},T]$ is the superposition of the individual percolations
\begin{equation} \label{csn}
\bfu_c(\bfx,t)=\sum_{k=0}^n \alpha^{(k)}\hbu^{(k)}(\bfx,t),\;\;\;\;\;p_c(\bfx,t)=\sum_{k=0}^n\hp^{(k)}(\bfx)/\beta^{(k)}(t),
\end{equation}
where $\hbu^{(0)}=\bbu$, $\alpha^{(0)}=1$; $\hp^{(0)}=\bp$, $\beta^{(0)}=1$.
 
The localisation and smoothing properties of the heat kernel (\ref{hhsp}) show that every $\hbu^{(k)}$ decays in time and in space, relative to its initial data $\hbu_0$. Moreover, all the percolations vanish on the plane $z=0$. Choose $\bfx=\bfx_s \in \Upo_h, t=T$ so that $\bbu(\bfx_s,T)\neq0$. We fix our $\alpha$'s by
\begin{equation*} 
\alpha^{(k)}\hbu^{(k)}(\bfx_s,T) = \bbu(\bfx_s,T)=A_s.
\end{equation*}
Thus, the velocity sum becomes $\bfu_c(\bfx_s,T)=(n'+1)A_s$, where $n'\leq n$, taking into account a finite number of (possible) zero $\hbu$.
By increasing the number of intervals over $[t^*,T]$ and by Cantor's diagonal process, we assert that
\begin{equation*} 
\bfu_c(\bfx_s,T)=\lim_{n\rightarrow \infty}(n'+1)A_s \; \longrightarrow \infty.
\end{equation*}
Furthermore, we choose a sub-interval, $[t_k,t_{k+1}]$, and divide it into $m$ equal parts, and repeat the whole analysis. The result is
\begin{equation*} 
\bfu'_c(\bfx,t)=\sum_{i=0}^m \alpha^{(i)}\hbu'^{(i)}(\bfx,t),\;\;\;t_k\leq t \leq t_{k+1},
\end{equation*}
which diverges at a point, ($\bfx_s,t_s$), as $m \rightarrow \infty$, as long as $\bbu(\bfx_s,t_s)$ does not vanish.

There is a considerable leeway to construct $\beta(t)$ in the sum (\ref{csn}), in order to deliver a divergent gradient in time, perhaps as an algebraic or exponential function. The perturbations given in (\ref{csn}) do not have to be singular, nevertheless, they do imply the existence of multiple backdrop flow-fields, starting from the identical initial conditions. 
Therefore, it makes no sense in talking about an {\it a priori} bound on the mean $(\bbu,\bp)$. The essence is that the Leray energy inequality (\ref{ienergy}) or similarity scalings are simply off-target, if not irrelevant. Without a background of the unpredictable percolated flows, claims of finite-time blow-up in the vicinity of the plane $z=0$ represent a misunderstanding of fluid motion. 

\section{Flow between plates}

We consider two steady planar flows between two parallel plates, namely plane Couette and Poiseuille flows. These flows have a long history of study, mainly in  stability theory. Without causing confusion, we will denote $\bfx_2{=}(x,y)$ and $\bfu{=}(u,v)$. The channel is unbounded: 
$\Upo_0: -\infty \:{<} \:x \:{<}\: \infty;\; 0\:{\leq}\: y\: {\leq}\: 1$.
The governing equation for real fluids $\mu>0$ has a simple form,
$\mu {\rd^2 u}/{\rd y^2} ={\rd p}/{\rd x}$. 
Plate boundary conditions yield a couple of basic profiles
\begin{align*} 
(\mbox{PCF})&: &u(x,0)&=0, & u(x,1)&=1; & \partial_x p&=0,   & u(x,y)&=y;\\
(\mbox{PPF})&: &u(x,0)&=0, & u(x,1)&=0; & \partial_x p&=-2\mu,&  u(x,y)&=y(1-y).
\end{align*}
We prevent ourselves from challenging the time-independence. To investigate the flow development with respect to perturbations, we consider the initial-boundary value problem of the full equations of motion (\ref{ns}). Specify the initial perturbation, $\hbu_0$, by
\begin{equation} \label{up}
\left( \begin{array}{c}
\hu_0  \\
\hv_0
\end{array} 
\right)(x,y)
	=\frac{1}{4 \pi} \int_{\Upo_0}  \log  \bigg( \:\frac{\cosh(x_-)-\cos(y_+)}{\cosh(x_-)-\cos(y_-)} \:\bigg) \left( \begin{array}{c}
\partial_{y'} \omega_0 \\
-\partial_{x'} \omega_0
\end{array} 
\right)	 \rd x' \rd y',
\end{equation}
where $x_-=\pi(x-x')$, $y_\pm=\pi(y \pm y')$, see, for instance, Carslaw \& Jaeger (1986). The properties of the logarithmic Green's function show $\hbu_0$ satisfies the no-slip on the plates. The smooth function, $\omega_0(x,y)$, is arbitrary. We require it to be any localised and bounded function around the origin. Since both the planar flows are time-independent, we denote the starting instant when the perturbation is first introduced by $t=0$. 

\subsection{Plane Couette and plane Poiseuille profiles}
%

%
Consider the perturbation on the linear profile, 
$(\bbu{+}\alpha \hbu, \hp/\beta)$, $\hbu{=}(\hu,\hv)$. In the governing equations (\ref{nss}), we replace $\bfu$ by $\bbu$, and set $\hbu$  the solutions of the planar diffusion equation, 
\begin{equation} \label{df2d}
	\partial_t \hbu - \nu \big(\partial_{xx} + \partial_{yy}\big) \hbu=0,
\end{equation}
subject to starting data $\hbu_0$, and homogeneous boundary value, $\hbu|_{\bdy_0}=0$. As the heat kernel takes the form, $H_2= H_1(x,x',t) {\times} H_f(y,y',t)$, where $H_f$ and $H_1$ are the kernels (\ref{hf}) and (\ref{h1}), we find the solution of (\ref{df2d})
\begin{equation} \label{spcf}
	\hbu(\bfx_2,t)= \int_{\Upo_0} \! H_2(\bfx_2,\bfy_2,t) \hbu_0(\bfy_2) \: \rd \bfy_2.
\end{equation}
The size of the percolation flow-field is much larger than unity for some values of $\alpha$, so that the rectilinear flow disappears soon after the perturbation is set-in. In two space dimensions, the right-hand function in (\ref{bsp}) takes a simpler form: 
\begin{equation} \label{bsp2}
\begin{split}
	\Delta \hp= {-} A_0 \: \Big(\: & \alpha \big(\: S_2[\hat{u}]\:\bfV_2[\hbu_x]+S_2[\hv]\:\bfV_2[\hbu_y]\: \big) \; + \\
	\quad & \hspace{7mm} 2 \big( \: S_2[u]\: \bfV_2[\hbu_x]+S_2[v]\: \bfV_2[\hbu_y]\: \big)\:\Big) \;\;(= -F_2),
\end{split}
\end{equation}
where, by analogy, $\bfV_2[\bfu_x]=(\partial_x u \;\partial_x v)'$, $
S_2[u]=(\partial_x u\; \partial_y u)$. Because of the linearity, the last term in (\ref{bsp2}) is reduced to $2 \partial_x \hv$. Given $\hbu$, the gradient counterpart of (\ref{bns}) actually upholds the momentum principle. In fact, the gradient is homogeneous on $\bdy_0$, as seen by inspecting (\ref{bns}). Pressure $\hp$ is evaluated from
\begin{equation} \label{ppcf}
	\int_{\Upo_0} N_2(\bfx_2,\bfy_2) F_2(\bfy_2, \hbu) \; \rd \bfy_2 + \hp_0,
\end{equation}
where, after an easy construction of the images, Neumann's function is 
\begin{equation} \label{n2}
	N_2(\bfx_2,\bfy_2)=-\frac{1}{4 \pi}\log  \Big(\big(\:\cosh(x_-)-\cos(y_-)\:\big)\:\big(\:\cosh(x_-)-\cos(y_+)\:\big) \Big).
\end{equation} 
%
%
In parallel to the linear profile, we seek a percolation of the form $\big(\bbu+\alpha \hbu,\;\bp+ \hp/\beta \big)$ on the parabolic mean flow. Given the class of initial data (\ref{up}), the expression (\ref{spcf}) remains unchanged for the description of the diffusive perturbations $\hbu(\bfx_2,t)$. Now, the forcing function for the pressure reads
\begin{equation*}
	F_2=- A_0\: \Big(\: \alpha \big(\: S_2[\hat{u}]\:\bfV_2[\hbu_x]+S_2[\hv]\:\bfV_2[\hbu_y]\: \big) +2(1 -2y) \partial_x \hv \: \Big),
\end{equation*}
while the pressure has the identical form of (\ref{ppcf}).

By concise computations, our analysis of the primitive Navier-Stokes dynamics shows that there exist no definitive disturbed states in each of these two shear flows, for any choice of $\omega_0$ of (\ref{up}). Once initiated, the basic profiles are completely outweighed by the imposing arbitrary perturbations which may be measured as infinitesimal or finite, albeit vaguely. The concept of disturbing established steady flows in the primitive setting is poised to be paradoxical.

\subsection*{Ill-posed initial-boundary value problem}

Instead of calling $\hu_0$ and $\hv_0$ of (\ref{up}) the perturbations, we set them as the initial data. Then one set of velocity solution $\hbu$ has the components
\begin{equation*}
\hu(x,y,t)=\hu_0(x,y)\hat{\alpha}(t),\;\;\;\hv(x,y,t)=\hv_0(x,y)\hat{\alpha}(t),
\end{equation*}
where $\hat{\alpha}$ is finite, regular but arbitrary, and $\hat{\alpha}(0){=}1$. Evidently, the flow-field remains solenoidal for all times $t > 0$ for every time-dependent function $\hat{\alpha}$. The following specifications of the pressure gradients ensure that the planar primitive equations are satisfied,
\begin{align*}
\partial_x \hp & = -\rho \big(\:\hu_0\hat{\alpha}' + \hat{\alpha}^2 (\hu_0 \partial_x \hu_0 +\hv_0 \partial_y \hu_0) - \hat{\alpha} \nu (\partial_{xx} + \partial_{yy}) \hu_0) \:\big),\\
\partial_y \hp & = -\rho \big(\: \hv_0\hat{\alpha}'  + \hat{\alpha}^2 (\hu_0 \partial_x \hv_0  +\hv_0 \partial_y \hv_0)  - \hat{\alpha} \nu (\partial_{xx} + \partial_{yy}) \hv_0) \:\big),
\end{align*}
for given $\omega_0$.

Our result does not contradict the global regularity of the two-dimensional equations established by Leray (1934$b$) because his proof of the well-posedness is based a maximum principle which bounds the scalar vorticity in $\real^2$. In this respect, the primitive formulation has been avoided.

\subsection{Spurious instability by normal-modes}

In the present section, we generalise the two basic uni-directional flows so that they depend on $x_3$ only, $\bar{\bfu}(x_3)$, and they are bounded between two infinite planes:
\begin{equation*}
\Upo_{\parallel}: \;\;\; -\infty < x_1,x_2 < +\infty, \;\;\; 0 \leq x_3 \leq 1.
\end{equation*}
By rewriting the basic profiles as $\bu=x_3$, and $\bu=x_3(1-x_3)$, we look for connections of the preceding analysis to the linear stability theory which postulates that any viscous flow may be decomposed into a basic profile plus perturbation
\begin{equation} \label{dcmp}
\bfu=\bbu+\bfu',\;\;\;p=\bp+p'.
\end{equation}
The normal-modes analysis in the instability theory is concerned with the excitation and attenuation of the wave-like perturbations
\begin{equation} \label{lst3}
	\bfu'(x_1,x_2,x_3,t)\propto \: \phi(x_3,t)\:\exp \big(\: \ri \tilde{\alpha} x_1 + \ri \tilde{\beta} x_2  \:\big),\\
\end{equation}
where $\tilde{\alpha}$ and $\tilde{\beta}$ stand for the respective wave numbers, and $\ri=\sqrt{{-}1}$. Only the real part of the complex representation in the eigen space is relevant to the motions. These wave-trains are known as instability waves which, presumably, exist and develop in space or in time, depending on whether the wave numbers, $\tilde{\alpha},\tilde{\beta}$, are complex or real.
 
In the domain of the parallel planes $\Upo_{\parallel}$, the complete Navier-Stokes equations (\ref{ns}) are invariant under the space translation in either the $x_1$-direction or the $x_2$-direction or both, namely
\begin{equation} \label{inv1}
	\bfu(x_1,x_2,x_3,t)=\bfu(x_1{-}a_1,x_2{-}a_2,x_3,t).
\end{equation}
It follows that the invariance principle must also apply to the combined flow, $\bar{\bfu}{+}\bfu'$. Because the basic flow itself is translation-invariant along the $x_1,x_2$ axes, the disturbances must preserve their analytical equivalence:
\begin{equation} \label{lstinv}
	A\: \phi(x_3,t)  \:\exp \big(\:\ri \tilde{\alpha} x_1 + \ri \tilde{\beta} x_2 \:\big) \equiv  
A \:\phi(x_3,t)\:\exp\big( \: \ri \tilde{\alpha} (x_1{-}a_1)+\ri \tilde{\beta} (x_2{-}a_2) \: \big),
\end{equation}
where $A$ is the proportional constant, and is independent of $x_1$ and $x_2$. For arbitrary choice of the translation scales $a_1$ or $a_2$, we must have 
\begin{equation} \label{nv}
\tilde{\alpha}=0,\;\;\; \mbox{and} \;\;\; \tilde{\beta}=0,
\end{equation}
for consistency of (\ref{lstinv}). At the planes, the boundary conditions such as the no-slip, are applied to the eigen part $\phi$. The non-existence deduction still holds regardless the presence of solid surfaces. As the perturbations are solenoidal, $\nabla{\cdot}\bfu'=0$, the wave-free constraints (\ref{nv}) show that $\partial \phi/\partial x_3=0$, which imposes strong restrictions on the eigen functions. Furthermore, the principle (\ref{inv1}) applies equally well to inviscid flows which are assumed to exist by formally setting $\nu=0$ in (\ref{ns}). Hence, the postulated sinuous instability waves (\ref{lst3}) are ill-defined in the framework of the exact incompressible Navier-Stokes-Euler dynamics, as their introduction implies incompatible physics.\footnote{The invariance law (\ref{inv}) also rules out the following mean profiles which depend linearly on the co-ordinates, 
\begin{equation*} 
\bu_i(x_j,t)=\sigma_{ij}(t)x_j+\bar{v}_i(t),\;\;\; x_j \in \real, \;\;\;(i,j=1,2,3).
\end{equation*}
A corresponding pressure can be found so that the equations of motion are satisfied. In $\real^3$, the proposed mean flows contain an infinite amount of energy for bounded $\sigma$'s. For arbitrary $\bfa$, the principle demands
\begin{equation*}
\bu_i(x_j,t)-\bu_i(x_j-a_j,t)=\sigma_{ij}(t)a_j=0,
\end{equation*}
which stipulates the discriminant: spatial strain rates $\sigma_{ij} \equiv 0$ at every time $t$. In rigorous fluid dynamics, we never recommend any exploration of spatially periodic perturbations on these profiles, as nothing insightful can be gained. In any event, these pathological flows would motivate false interpretations on fluid flows, such as the `generalised Kelvin modes'. There have been some confusions about Kelvin's paper (1887) in technical literature. Kelvin had worked exclusively on bounded uni-directional plane Couette flow between $0 \leq y \leq b$, though he mentioned an extension, $b \rightarrow \infty$, on p.$191$ of his paper. Apparently, it was Hopf (1952) who first considered the boundary-free turbulence (see his \S 7) as a limit of spatially periodic flow. However, he expressed strong reservations about the unbounded behaviours at large distances, when the period tends to infinity. Given our current knowledge, we assert that, as a result of careful analyses, there are no infinite-energy viscous motions. Consequently, the invariance principle is an equivalent statement to the non-existence of the Navier-Stokes solutions in periodic domains.}

\subsection{Unpredictability due to disturbances}

By the decompositions (\ref{dcmp}), the superimposed flow of a basic profile and a perturbation satisfies the equations of motion
\begin{equation} \label{nsprt}
\partial_t \bfu' + (\bfu'{\cdot} \nabla) \bfu'+ (\bbu{\cdot}\nabla)\bfu' + (\bfu'{\cdot}\nabla)\bbu = - \nabla p' /\rho+ \nu \Delta \bfu',\;\;\;\nabla{\cdot}\bfu'=0,
\end{equation}
because the basic flow, $(\bbu,\bp)(\bfx,t)$, is a solution. As in many practical applications, it has been claimed that equations (\ref{nsprt}) can be solved by numerical methods, for given initial data. Let us call these (primitive) solutions $(\tbu,\tp)$, and they are valid over time interval $[t_0,T^*]$. 

Choose a sub-interval $[t_0,T]$, $T<T^*$. Divide the time interval $[t_0,T]$ into $n$ equal parts: $\Delta t = (T-t_0)/n$ or $t_k= k\Delta t, \;k=1,2,\cdots,n$. We consider the first percolation on the mean flow field at $t_0$,
\begin{equation} \label{pc0}
(\bbu+\alpha^{(1)}\hbu^{(1)},\; \bp+\hp^{(1)}/\beta^{(1)}).
\end{equation}
Let the operator ${\bf H}(\bfu)=(\partial_t  -  \nu \Delta )\bfu$.
Then the percolation is determined by solving the heat equation with a buffer source
\begin{equation} \label{pc1}
{\bf H}(\hbu^{(1)})= \bfb,
\end{equation}
subject to (arbitrary) solenoidal initial data, $\hbu_p$. Because of the sandwiched geometry, there are useful symmetries in the Green's function $H_3$ which is governed by 
\begin{equation*}
{\bf H}(H_3)=-\delta(x-x')\delta(y-y')\delta(z-z')\delta(t).
\end{equation*}
Taking into account the non-slip on $\partial \Upo_{\parallel}$, we find that $H_3$ has a composite form
\begin{equation*}
H_3(\bfx, \bfx',t)= H_1(x,x',t) \times H_1(y,y',t) \times H_f(z,z',t).
\end{equation*}
Thus the solenoidal diffusion is guaranteed,
\begin{equation} \label{d1}
\hbu^{(1)}(\bfx,t)=\int_{\Upo_{\parallel}} \!\! H_3 (\bfx{-}\bfx',t) \hbu_p(\bfx') \rd \bfx'+\int_{t_0}^{t_1}\!\!\!\int_{\Upo_{\parallel}} \!\! H_3 (\bfx{-}\bfx',t{-}s) \bfb(\bfx') \rd \bfx' \rd s,
\end{equation}
where $\hbu_p \in L^1_{loc}(\Upo_{\parallel})$, for any choice of the source $\bfb$ within the constraints imposed on the buffer functions (\ref{buff}). One of the important properties of the kernel $H_3$ is that it is capable of mollifying a large class of the initial data. In experimental work, various methods have been developed to disturb the mean flow, including wave generators, acoustic excitations by loud speakers, flow impulses, controlled directional jets, and others. All these artificial excitations are well-covered in the theory of $L^1_{loc}$ space. Computationally, the diffusive component $\hbu^{(1)}$ can be recognised as the numerical errors from discretisation.

The first equality of (\ref{nsprt}) is fully maintained if
\begin{equation} \label{ppc0}
\nabla \hp^{(1)} = - {\bf q}^{(1)}(\hbu^{(1)},\bbu),
\end{equation}
where the non-linear terms are grouped into 
\begin{equation*}
{\bf q}^{(i)}(\bfu,\bfv) = \rho \alpha^{(i)} \beta^{(i)} \big(\: \alpha^{(i)} (\bfu{\cdot} \nabla) \bfu+ (\bfv{\cdot}\nabla)\bfu + (\bfu{\cdot}\nabla)\bfv \:\big) + \rho \beta^{(i)} \bfb.
\end{equation*}
The actual $\hp$ can be found by solving the problem
$\Delta \hp^{(i)} = - (\nabla{\cdot}{\bf q}^{(i)})$,
subject to inhomogeneous Neumann data of (\ref{ppc0}). Some algebra gives
\begin{equation*} 
\begin{split}
\hp^{(i)}(\bfx)&=\int_{\Upo_{\parallel}} \!\! N_3(\bfx,\bfx') (\nabla{\cdot}{\bf q}^{(i)})(\bfx')\: \rd \bfx'  \\
\quad & \;\;+ \rho \beta^{(i)} \! \int_{\real^2} \! \Big( N_3(\bfx,\bfx_2',0)\bfb(\bfx_2',0) -N_3(\bfx,\bfx_2',1)\bfb(\bfx_2',1) \Big) \rd \bfx_2'+ \hp_0,
\end{split}
\end{equation*}
where
\begin{equation*}
N_3=\frac{1}{4 \pi}\sum_{k=-\infty}^{\infty}\Big(\: \frac{1}{R_+}+\frac{1}{R_-}\:\Big)
\end{equation*}
with the image locations at $R^2_{\pm}= (x{-}x')^2+(y{-}y')^2+(z{\pm} z' {-} 2k)^2$.
Hence our first companion solution deviates the mean motion
\begin{equation*}
(\bfu_c,p_c)(\bfx,t)=\Big(\: \bbu+\alpha^{(1)}\hbu^{(1)},\; \bp+\hp^{(1)}/\beta^{(1)}\:\Big)(\bfx,t).
\end{equation*}
In general, the postulated solution $(\tbu,\tp)$ and the companion do not coincide over $[t_0,t_1]$, even though both are perturbed by the same initial data $\hbu_p$.

Next we consider a series of consecutive percolations
\begin{equation*}
\Big(\alpha^{(k)}\hbu^{(k)},\; \hp^{(k)}/\beta^{(k)}\Big),\;\;\;k=2,3,\cdots,
\end{equation*}
on the companion flow
\begin{equation*}
(\bfu^{(k-1)}_c,p^{(k-1)}_c)=\bigg(\bbu+\sum_{i=1}^{k-1}\alpha^{(i)}\hbu^{(i)},\; \bp+\sum_{i=1}^{k-1}\hp^{(i)}/\beta^{(i)}\bigg),
\end{equation*}
over $t_{k-1} \leq t \leq t_{k}$.
The $k$-th percolation is found from ${\bf H}(\hbu^{(k)})=\bfb$, subject to data $\bfu_c^{(k-1)}$, and it has the form (\ref{d1}). By the same token, we must have 
\begin{equation*}
\nabla \hp^{(k)} = - {\bf q}^{(k)}\big( \hbu^{(k)},\bfu^{(k-1)}_c \big).
\end{equation*}
Over the last interval $[t_{n-1},t_n]$, the basic profile $(\bbu,\bp)$ has an arbitrary companion
\begin{equation*} 
(\bfu^{(n)}_c,\nabla p^{(n)}_c)(\bfx,t)=\bigg(\bbu+\sum_{k=1}^{n}\alpha^{(k)}\hbu^{(k)},\; \nabla \bp+\sum_{k=1}^{n}\nabla \hp^{(k)}\bigg)(\bfx,t).
\end{equation*}
By increasing the number of sub-intervals, $n \rightarrow \infty$, we envisage an experimental condition of continuously introducing perturbation $\hbu_p$ into the mean flow. A typical example of this situation is found in all wind tunnels where the free-stream contains high levels of turbulence or acoustic noise. Another example includes the macroscopic realisation of thermal fluctuations or Brownian motion which is an inherent part in all flow states. Evidently, the claimed primitive solution, $(\tbu,\nabla\tp)$, can never be unique over $[t_0,T]$, in spite of the origin of the perturbations. By deduction, the primitive theory implies that any mean flow will be strongly distorted immediately after the start. This theoretical account does not reconcile with experimental facts. 

Our evaluations so far have been carried out by exact computations. We have kept all the terms in equation (\ref{nsprt}). In linear stability theory, it has been assumed that the term, $(\bfu'. \nabla) \bfu'$, is small compared to other terms in (\ref{nsprt}). Neglecting this non-linear term, we have 
\begin{equation} \label{lns}
\partial_t \bfu' + (\bbu{\cdot}\nabla)\bfu' + (\bfu'{\cdot}\nabla)\bbu = - \nabla p'+ \nu \Delta \bfu',\;\;\;\nabla{\cdot}\bfu'=0.
\end{equation}
In principle, the linearisation is made at the expense of the momentum conservation; there are no justifications for ignoring the non-linearity unless one tolerates absurd physics. Nevertheless, let the solution of (\ref{lns}) be $(\bbu_L,\bp_L)$ over its validity interval $[t_0,T]$. This basic flow may have been obtained, for instance, by a normal-modes analysis of linear instability or a procedure of asymptotic expansion at high Reynolds numbers. 
Starting with the first percolation (\ref{pc0}), ${\bf H}(\hbu_L^{(1)})=\bfb$ remains unchanged. Clearly, the linearisation does not modify our diffusion problem for $\hbu_L^{(k)}$. Only minor adjustments are required for the compensating pressure $\hp_L$. We need
\begin{equation*}
{\bf q'}^{(i)}(\bfu,\bfv) = \rho \alpha^{(i)} \beta^{(i)} \big(\: (\bfv{\cdot}\nabla)\bfu + (\bfu{\cdot}\nabla)\bfv \:\big)+\rho \beta^{(i)} \bfb,
\end{equation*}
instead of ${\bf q}$. Without repeating the whole analysis for consecutive percolations, we deduce that the `linearised' companion to $(\bbu_L,\bp_L)$ contains the similar arbitrariness
\begin{equation*} 
(\bfv^{(n)},\nabla q^{(n)})(\bfx,t)=\bigg(\bbu_L+\sum_{k=1}^{n}\alpha^{(k)}\hbu_L^{(k)},\; \nabla \bp_L+\sum_{k=1}^{n}\nabla \hp^{(k)}_L\bigg)(\bfx,t),
\end{equation*}
where
$\nabla \hp^{(k)}_L = - {\bf q'}^{(k)}\big(\hbu_L^{(k)},\bfv^{(k-1)} \big)$.
As in the case of the upper half space, infinitely many consecutive percolations are able to saturate any basic profile with the aid of selective $\alpha$'s and $\beta$'s. What we are trying to emphasise is the fact that, the connection of the sinusoids (\ref{lst3}) to fluid motions are only found in the linearised theory. The key issue is that these wave-like instabilities comprise just a small subset of the allowable $L^1_{loc}$ perturbations which take any physically-relevant forms, for instance, spatiotemporal aperiodic fluctuations. To justify specific normal modes by evoking Fourier transforms is a half-done argument. Indeed counter examples to normal-modes analyses are well-known, as there are flows whose spectra are ostensibly incomplete. In such circumstances, the expansion theory of instability modes is at least inapplicable, as the continuum flow scales do not amend to wave descriptions.

It is an inconvenient reality that our description of fluid motions solely in terms the primitive variables precipitates the inevitability of the momentum partition, which gives rise to successive percolations under the auspices of the buffer $\bfb$. The grave consequence is that any disturbed state, $(\bbu+\bfu',\;\bp+p')$, amounts to indeterminacy. Hence, those methods of solving the primitive (\ref{nsprt}), such as linear or weakly non-linear instability analyses, non-modal optimised energy approaches, Fourier or eigen-function expansions, high Reynolds-number asymptotic approximations, can hardly be promising.

\subsection{Boundary layers}

A thin viscous layer develops on the surface of a plate when it is suddenly set into motion. Theoretically, the thickness of the layer can be estimated from the impulsively-started Stokes' flow on an infinite plane, 
\begin{equation*}
\partial_t u-\nu \partial_{yy} u=0,
\end{equation*}
subject to $u_0(y,0){=}0, u(0,t\geq0){=}U_{\infty}$, and $u(\infty,t){=}0$. The solution of this equation is given by
\begin{equation*}
u(y,t)=U_{\infty}\int_0^t \!\! T(t{-}s)\rd s = \frac{y U_{\infty}}{\sqrt{4 \pi \nu}} \int_0^t \frac{\exp\big({-}{y^2}/{\big(4\nu(t{-}s)\big)}\big) }{(t{-}s)^{3/2}}\rd s =\frac{y U_{\infty}}{\sqrt{4 \pi \nu}} I(t),
\end{equation*}
where, from (\ref{h1}),
\begin{equation*}
T(t-s)={\partial}( H_1(y,y',t{-}s) - H_1(y,{-}y',t{-}s))/{\partial y'}\big|_{y'=0}.
\end{equation*}
The integral, $I(t)$, is a measure of time scale. A fluid particle on the plate impulsively moves with the plate to a location $L$ in time $L/U_{\infty}=t_L$. In fluids $\nu \ll 1$, the exponential function is $\sim O(1)$ for $y^2/(\nu t) \ll 1$. Asymptotically, $t_L$ is the characteristic time 
\begin{equation*}
I(t) \sim \int_0^t (t-s)^{-3/2} \rd s \approx t^{-3/2} \int_0^{t} \big(1 - {3s}/{(2t)}\big) \rd s \approx t_L^{-1/2} = \sqrt{{U_{\infty}}/{L}}.
\end{equation*}
Close to the plate, $u/U_{\infty} \sim O(1)$, the wall layer, $\delta \sim y$, is called a boundary layer. For an observer moving with the plate, the thin viscous shear  
\begin{equation} \label{bldel}
\delta  \sim O(\sqrt{\nu}) \;\;\;\;\;\mbox{or}\;\;\;\;\; {\delta}/{L} \sim \sqrt{{\nu}/{L U_{\infty}}} = Re^{-1/2}_L,
\end{equation}
for Reynolds' number $Re_L \gg 1$.

The boundary layer theory of Prandtl (1904) is an approximation of the planar Navier-Stokes dynamics. For a steady boundary layer over a semi-infinite flat-plate, the governing equations are 
\begin{equation} \label{bl0}
\frac{\partial u}{\partial t}+u\frac{\partial u}{\partial x} + v\frac{\partial u}{\partial y} = - \frac{1}{\rho}\frac{\partial p}{\partial x} + \nu \frac{\partial^2 u}{\partial y^2},\;\;\;\;\;\; \frac{\partial u}{\partial x} +\frac{\partial v}{\partial y}=0.
\end{equation}
The boundary conditions are $
u(x,0){=}v(x,0){=}0; u {\rightarrow} U_{\infty}, v {\rightarrow} 0,y {\rightarrow} \infty$,
and upstream compatibility $u(0,y){=}U_{\infty}, v(0,y){=}0$. In the scaling $U^2_{\infty}/L$, $v$-momentum is of small order, and $\partial p/\partial y \sim O(\delta/L)$. Since there are three unknowns in (\ref{bl0}), a supplementary condition must be imposed. It is argued that, at the edge of the layer, velocity $u$ presumes the value of the inviscid free-stream $U_{\infty}$. Thus the stream-wise pressure satisfies
\begin{equation} \label{blp} 
{\rd p}/{\rd x} = - \rho U_{\infty}\:{\rd U_{\infty}}/{\rd x}.
\end{equation}
Effectively, the gradient in (\ref{bl0}) is fixed for every given $U_{\infty}{=}U_{\infty}(x)$. The well-studied case is a Blasius profile on a flat-plate, $U_{\infty}{=}U_0{=}\const$ In steady flows, a number of free-streams have been proposed, for example, $U_{\infty}(x) {=} U_0 \big(x/L\big)^m$ for an integer $m$. In the flow over a circular cylinder, the initial flow and free-stream are matched, $u_0(x,0)=U_{\infty}(x)=U_0\sin(x/L)$, and hence viscous layers develop in the vicinity of the forward stagnation point. The expectation is that, from the solutions of (\ref{bl0}), boundary layers in favourable and adverse gradients may be elucidated. However, the approximations (\ref{bl0}) are a reductionist's view of fluid motion and are not always reasonable. For instance, elevating passive pressure (\ref{ppoi}) to regulative control (\ref{blp}) introduces a misconception that an adverse pressure gradient causes a boundary layer to separate.

Incidentally, the supposition of a fixed gradient circumvents the non-uniqueness of the primitive formulation, as it closes the partition loop of the momentum. But, considerations of perturbing a Blasius flow $\bbu_{BL}$, in the sense of (\ref{dcmp}), necessarily involve an analogous analysis to (\ref{nsprt}) or (\ref{lns}). In practice, a typical flat-plate boundary layer has a stream-wise length $L \sim O(1)$ behind the leading edge at moderate Reynolds numbers (cf. the estimate (\ref{bldel})). Comparisons between measurements and linearised calculations may seem to be impressive. Even so, we must treat any of the apparent successes with caution, as these limited test realisations by no means elude the analytical indeterminacy. Accordingly, the contrast between convective and absolute instabilities is inconsequential.\footnote{Our findings of the present section are consistent with experiment and careful computations. The classic experiment on the boundary-layer transition by Schubauer \& Skramstad (1947) conclusively demonstrated that perturbations have no correlations with the measured transition location. Instead of the primitive variables, solutions must be sought in the continuum theory of the full vorticity dynamics. In applications, we do not regard unqualified ripples and stitches in vorticity iso-contours as evidence of the flow instabilities. Yet the critical layers are ill-defined and hence dissolved. Highly skewed vortices haul retarded flow elements away from solid surfaces, leading to separations. Convergence in discretisation and successive asymptotes is of importance. We exemplify these issues by a dipole-wall simulation, see figures $7$ to $11$ of Lam (2018$a$).}

\section{Cylindrical polar co-ordinates}

Denote the velocity components by $\bfu=(u_r,u_{\theta},u_z)$ in the cylindrical co-ordinates $\bfz=(r,\theta,z)$. The primitive equations of motion are
\begin{equation} \label{nsc}
\begin{split}
		\frac{\partial u_r}{\partial t}+ (\bfu{\cdot}\nabla)u_r - \frac{u^2_{\theta}}{r} & =-\frac{1}{\rho}\frac{\partial p}{\partial r}+ \nu \bigg(\Delta_1 u_r -\frac{2}{r^2}\frac{\partial u_{\theta}}{\partial \theta}\bigg),\\
		\frac{\partial u_{\theta}}{\partial t} + (\bfu{\cdot}\nabla)u_{\theta} + \frac{u_r\: u_{\theta}}{r}& =-\frac{1}{\rho \: r}\frac{\partial p}{\partial \theta} + \nu \bigg( \Delta_1 u_{\theta} + \frac{2}{r^2}\frac{\partial u_r}{\partial \theta}\bigg),\\
		\frac{\partial u_z}{\partial t }+ (\bfu{\cdot}\nabla)u_z &=-\frac{1}{\rho}\frac{\partial p}{\partial z} +  \nu \:\Delta_0 u_z, 
		\end{split}
\end{equation}
and the continuity,
\begin{equation} \label{divc}
		\nabla'{\cdot} \bfu = \frac{\partial u_r}{\partial r}+\frac{u_r}{r}+\frac{1}{r}\frac{\partial u_{\theta}}{\partial \theta} + \frac{\partial u_z}{\partial z}=0,
\end{equation}
where the operators are defined by $\Delta_1=\Delta_0-1/r^2$,  
\begin{equation*}
	\Delta_0 = \frac{\partial^2}{\partial r^2}+ \frac{1}{r} \frac{\partial }{\partial r} + \frac{1}{r^2}\frac{\partial^2}{\partial \theta^2}+ \frac{\partial^2 }{\partial z^2}.
\end{equation*}
The gradient operator, $
\nabla = {\partial}/{\partial r }+{\partial}/{(r\partial \theta)}+{\partial}/{\partial z}$.
Equations (\ref{nsc}) and (\ref{divc}) are made dimensionless by the SI units. To derive the analogous expression of (\ref{ppoi}) for the pressure, we take into account the fact that the time-dependent terms must vanish in view of the incompressibility hypothesis. In practice, we look for an expression containing $\Delta_0p $ or $\Delta_1 p$. The following identities facilitate the subsequent algebraic simplifications:
\begin{equation} \label{s3}
\begin{split}
\frac{\partial}{\partial r }\big(\Delta_0 u_r\big)&= \Delta_0\Big( \frac{\partial u_r}{\partial r }\Big)-\frac{1}{r^2}\frac{\partial u_r}{\partial r }- \frac{2}{r^3 }\frac{\partial^2 u_r }{\partial \theta^2 },\\
\frac{1}{r }\big(\Delta_0 u_r\big)&= \Delta_0\Big( \frac{u_r}{r}\Big)+\frac{2}{r^2 }\frac{\partial u_r}{\partial r }-\frac{u_r}{r^3 },\\
\frac{1}{r}\frac{\partial}{\partial \theta }\big(\Delta_0 u_{\theta}\big)&= \Delta_0\Big(\frac{1}{r} \frac{\partial u_{\theta}}{\partial \theta }\Big)+\frac{2}{r^2}\frac{\partial^2 u_{\theta}}{\partial r \partial \theta}-\frac{1}{r^3}\frac{\partial u_{\theta}}{\partial \theta },
\end{split}
\end{equation}
where each term of the left-hand side can be found by manipulating the governing equations. By direct computations, the momentum equations in (\ref{nsc}) then show that the diffusion-related terms drop out due to the continuity.
Effectively, we need to simplify the expression
\begin{equation*}
\Big(\frac{\partial}{\partial r}+\frac{1}{r}\Big)\Big((\bfu{\cdot}\nabla)u_r-\frac{u^2_{\theta}}{r} \Big) + \frac{1}{r}\frac{\partial}{\partial \theta}\Big((\bfu{\cdot}\nabla)u_{\theta}+\frac{u_r\:u_{\theta}}{r}\Big)+ \frac{\partial}{\partial z}\Big((\bfu{\cdot}\nabla)u_z\Big),
\end{equation*}
which equals to $-\Delta_0 p/\rho$.
It is found that the lengthy identity can be paraphrased by subtracting a variant of the continuity, $(\bfu.\nabla)(\nabla'{\cdot}\bfu)$,
so that all the terms containing second derivatives cancel out.
After further reductions and some rearrangements, we obtain the Poisson equation for the pressure
\begin{equation} \label{ppoic}
\begin{split}
- \frac{1}{\rho}\Delta_0 p & = \Big(\frac{\partial u_r}{\partial r}\Big)^2+\Big(\:\frac{1}{r}\frac{\partial u_{\theta}}{\partial \theta}+\frac{u_r}{r}\Big)^2+\Big(\frac{\partial u_z}{\partial z}\Big)^2\\
\quad & \hspace{10mm} + 2 \Big(\: \frac{1}{r}\frac{\partial u_r}{\partial \theta}\frac{\partial u_{\theta}}{\partial r} + \frac{1}{r}\frac{\partial u_{\theta}}{\partial z}\frac{\partial u_z}{\partial \theta} + \frac{\partial u_r}{\partial z}\frac{\partial u_z}{\partial r}- \frac{u_{\theta}}{r}\frac{\partial u_{\theta}}{\partial r}\:\Big)=P(\bfu).
\end{split}
\end{equation}

\subsection{Non-existence of Burgers' vortex}

Flows in the domain,
\begin{equation*}
 \Upo_{z}: \;\;\; 0 \leq r < \infty; \;\;\; 0 \leq \theta < 2 \pi;\;\;\; - \infty < z < +\infty,
\end{equation*}
must obey the invariance principle (\ref{inv}) in the $z$-direction. Burgers (1948) suggested a vortex motion whose velocities are 
\begin{equation} \label{bvx}
u_r=- \gamma r;\;\;\;u_{\theta}=\phi(r,t);\;\;\;u_z= 2 \gamma z,
\end{equation}
where $\gamma$ denotes a constant, representing a strain rate as it has the dimension $[\texttt{s}^{-1}]$. The invariance law in the $z$-direction is expressed in the identity,
\begin{equation*}
u_z(z)-u_z(z-a_3)= 2 \gamma a_3 =0,
\end{equation*}
for arbitrary translation $a_3$. Thus, Burgers' vortex cannot exist for non-zero $\gamma$. Few variants of (\ref{bvx}), for instance, $\gamma=\gamma(t)$ or $\gamma=\gamma f(r)$, have been proposed. These revised vortex flows set the scene for a number of strained vorticity fields, which were treated as building elements in simulating turbulence. By inference, the conclusions of these works have marginal validity, and large parts of them, if not all, must be reappraised for important applications. In finite energy flows, the rule is that the Navier-Stokes equations do not admit solutions of power series 
\begin{equation*}
u_i \propto (x_i-a_i)^k,\;\;\;(i=1,2,3)
\end{equation*}
in homogeneous direction $x_i$ for non-zero integer $k$.

\subsection{Rotational normal-modes for Hagen-Poiseuille flow}

In the infinite circular pipe of unit radius,
$\Upo_1: 0 \leq r \leq 1; \;0 \leq \theta < 2 \pi;\; z \in \real$,
we consider a uni-directional flow $\bu_z=\bu_z(r)$ with zero radial and circumstantial pressure gradients. The continuity is satisfied automatically, and the $z$-momentum of (\ref{nsc}) reduces to $\mu (\partial_{rr} +\partial_r/r)\bu_z = \partial \bar{q}/\partial z $. Applying the no-slip on wall $r=1$, we obtain the steady velocity, 
\begin{equation} \label{pipe}
	\bar{\bfu}=\big(\bu_r,\; \bu_{\theta}, \; \bu_z
	\big)=\big(\:0,\;\;0,\;\;1-r^2\:\big),
\end{equation}
which is driven by an axial pressure gradient, $\partial_z \bar{q}=-4 \mu$. In laboratory experiments, the steady-state profile (\ref{pipe}) may be realised at locations a few diameters downstream of a streamlined entry, as long as the transient of starting conditions becomes settled. In this way, the comparison with the theory tentatively removes the dependence of initial conditions. Given the subtleties and difficulties in fluid experiments, this unitary case should never be generalised into a conception that initial conditions are of secondary importance. Practically, we may devise many entry conditions, and hence downstream flows are strong functions of the applied external pressure gradient. 

In theory, parabolic distribution (\ref{pipe}) is applicable over the entire domain $\Upo_1$ and is thus treated as a mean flow. In parallel to the wave-form (\ref{lst3}), stability analysis considers possible amplifying effects of the disturbances 
\begin{equation} \label{dst3}
	(u'_r,u'_{\theta},u'_z,p')(r,\theta,z,t)\propto(\tilde{u}_r,\tilde{u}_{\theta},\tilde{u}_z,\tilde{p})(r,t)\: \exp\big(\:\ri\:(\:\tilde{\gamma} z + n \theta\:) \:\big),
\end{equation}
where $n=0,\pm1,\pm2,\cdots$, so as to ensure rotational periodicity. The principle of translation invariance asserts that the dynamic compatibility must be observed at all times, or the following expression must make sense
\begin{equation*}
\exp\big(\:\ri\:(\:\tilde{\gamma} z + n \theta\:) \:\big) 
= \exp\big(\:\ri\:(\:\tilde{\gamma} (z-a_3) + n \theta\:) \:\big).
\end{equation*}
In the exact non-linear Navier-Stokes dynamics, the travelling waves (\ref{dst3}) cannot exist for arbitrary $a_3$ unless $\tilde{\gamma}=0$. As a result, the expected perturbed structures are less interesting, as they degenerate into fixed periodic functions in $\theta$. 

By integration, the actual pressure driving the flow (\ref{pipe}) is found to have an analytic expression $\bar{q}=-4 \mu z + \bar{q}_0(t)$. From a strictly analytical point of view, we may articulate that there can be no pressure of this form which satisfies the invariance law. However, we shall ignore this indeterminate quantity, and stick to the fact that the momentum conservation involves only the gradient.

\subsection{Ambiguity in pipe flow}

Let us introduce an add-on flow, $\hbu$, superimposed on flow (\ref{pipe}), 
\begin{equation} \label{ap}
	(\hbu,\hat{q})(r,\theta,z,t)=(\hu_r,\hu_{\theta},\hu_z,\hat{q})(\bfz,t),
\end{equation}
where the velocity vanishes at the pipe wall, $\hbu(1,\theta,z,t)=0$, and decays to zero as $z \rightarrow \pm \infty$. We substitute the sum, 
$\bar{\bfu}+\alpha\hbu$, $\bar{q}+\hat{q}/\beta$,
into equations (\ref{nsc}) and (\ref{divc}). Since $(\bbu,\baq)$ is a solution pair, the analogous system to (\ref{nss}) reads
\begin{equation} \label{nsd}
\begin{split}
 \partial_t \hu_r + f_r &=-\big({\partial \hq}/{\partial r}\big)/A_0+ \nu \big(\Delta_1 \hu_r -{2}{r^{-2}}{\partial \hu_{\theta}}/{\partial \theta}\big),\\
 \partial_t \hu_{\theta} + f_{\theta} &=-\big({\partial \hq}{\partial \theta}\big)/(r A_0) + \nu \big( \Delta_1 \hu_{\theta} + {2}{r^{-2}}{\partial \hu_r}/{\partial \theta}\big),\\
 \partial_t \hu_z + f_z &=-\big({\partial \hq}/{\partial z}\big)/A_0 + \nu \:\Delta_0 \hu_z, \\
\nabla'.\hbu&=0.
\end{split}
\end{equation}
where the non-linear terms are grouped into
\begin{align*}
 f_r & =\alpha (\hbu{\cdot}\nabla)\hu_r + (1-r^2) \partial \hu_r/\partial z - \hu^2_{\theta}/r,\\ 
 f_{\theta} & = \alpha(\hbu{\cdot}\nabla)\hu_{\theta} + (1-r^2) \partial \hu_{\theta}/\partial z + \hu_r \hu_{\theta}/r, \\
  f_z & = \alpha(\hbu{\cdot}\nabla)\hu_z + (1-r^2) \partial \hu_z/\partial z - 2r \hu_r.
\end{align*}
The Neumann data can be inferred from (\ref{nsd}), the equation gives the pressure 
\begin{equation} \label{ppcw}
\Delta_0\hat{q} = -A_0 \big({\partial f_r}/{\partial r}+{f_r}/{r} + {\partial f_{\theta}}/{\partial \theta}/r+ {\partial f_z}/{\partial z}\big)=-\hat{Q}_p(\bbu, \hbu,\bfz).
\end{equation}

Choose any buffer $\bfb=(b_r,b_{\theta},b_z)(\bfz,t)$, as postulated. The components of the percolation are diffusions
\begin{equation} \label{apd}
{\bf h}_1(\hu_r)=-{2}{r^{-2}}{\partial \hu_{\theta}}/{\partial \theta}+b_r,\;\; {\bf h}_1(\hu_{\theta})={2}{r^{-2}}{\partial \hu_r}/{\partial \theta}+b_{\theta},\;\;{\bf h}_0(\hu_z)=b_z,
\end{equation}
where ${\bf h}_0{=}\partial_t {-} \nu \Delta_0$, and ${\bf h}_1{=}{\bf h}_0{+}1/r^2$.
Because the Hagen-Poiseuille flow is steady, we are free to set the switch-on time of percolations, say, $t=0$. The initial data are
\begin{equation} \label{apic}
	\big(\:\hu_r,\hu_{\theta},\hu_z\:\big)(\bfz,t=0)=\big(\:\hu_0,\;\hv_0,\;\hw_0\:\big)(\bfz).
\end{equation}
The solution of the third equation in (\ref{apd}) is expressed as
\begin{equation} \label{uz}
	\hu_z(\bfz,t) =\int_{\Upo_1}\!\!K_0(\bfz,\bfz',t) \: r' \hw_0(\bfz') \: \rd \bfz'+ \int_0^t \!\!\int_{\Upo_1} \!\!K_0(\bfz,\bfz',t{-}s) \: r'b_z(\bfz',s) \:\rd \bfz' \rd s,
\end{equation}
where $K_0(\bfz,\bfz',t)=H_1(z,z',t)G_0(r,\theta,r',\theta',t)$, with the cross-section Green function being
\begin{equation} \label{g0}
 \begin{split}
	G_0& =\frac{1}{\pi}\: \sum_{i=1}^{\infty} \frac{1}{(J_1(\kappa_{0i}))^2} \: J_0(\kappa_{0i} r) \: J_0(\kappa_{0i} r') \: \exp(-\kappa^2_{0i} \nu t) \\
	\quad & \hspace{10mm} + \frac{2}{\pi} \sum_{n,i=1}^{\infty}  \frac{\cos(n(\theta-\theta'))}{(J_{n+1}(\kappa_{ni}))^2}\:J_n(\kappa_{ni} r)\:J_n(\kappa_{ni} r')\: \exp(-\kappa^2_{ni} \nu t),
	\end{split}
\end{equation}
where $J_n$ stands for the Bessel function of order $n=0,1,2,\cdots$. Each of the Bessel functions, $J_n(x)$, has infinitely many zeros in the interval $x>0$. In every summation over $n$, $\kappa_{ni}$ is the $i$th positive root of the function $J_n(\kappa)=0$. 

Let $u_{\pm}=\hu_r \pm \hu_{\theta}$. The first two equations of (\ref{apd}) are reduced to the system
\begin{equation} \label{upm}
\left.
\begin{aligned}
	(\partial_t & - \nu \Delta_1 ) u_{+}= 2\nu\:\partial_{\theta} u_{-}/r^2+b_r +b_{\theta},\\
	(\partial_t & - \nu \Delta_1 ) u_{-}= -2\nu \: \partial_{\theta} u_{+}/r^2+b_r - b_{\theta}.
	\end{aligned}
\hspace{5mm}\right\}
\end{equation}
Since the diffusion operator on the left is separable, we readily establish the analogous kernel, $K_1(\bfz,\bfz',t)=H_1(z,z',t) G_1(r,\theta,r',\theta',t)$, where the Green function $G_1$ is
\begin{equation*} 
	G_1(r,\theta,r',\theta',t)=\frac{2}{\pi}\sum_{n,i=1}^{\infty}  \frac{\cos(n(\theta-\theta')) }{(J_{n+1}(\kappa_{ni}))^2} \: J_{n}(\kappa_{ni} r) \: J_{n}(\kappa_{ni} r') \: \exp(-\kappa^2_{ni} \nu t).
\end{equation*}
Transforming (\ref{upm}), the solutions are defined in the coupled integral equations:
\begin{equation} \label{ie2}
\left.
\begin{aligned}
u_+(\bfz,t) &+ 2 \nu \int_0^t\!\!\int_{\Upo_1} K_2(\bfz,\bfz',t{-}s)/r'\: u_-(\bfz',t{-}s)\:\rd \bfz' \rd s = F_+(\bfz,t), \\
u_-(\bfz,t) &- 2 \nu \int_0^t\!\!\int_{\Upo_1} K_2(\bfz,\bfz',t{-}s)/r'\: u_+(\bfz',t{-}s)\:\rd \bfz' \rd s = F_-(\bfz,t), \\
\end{aligned}
\hspace{5mm}\right\}
\end{equation}
where the integral terms are obtained by integration by parts, taking into account the no-slip condition. Kernel $K_2(\bfz,\bfz',t)=H_1(z,z',t) G_2(r,\theta,r',\theta',t)$. The drivers on the right, $F_{\pm}$, are known
\begin{equation*} 
\int_{\Upo_1}\!\! K_1 (\bfz,\bfz',t)r' \: (\hu_0 \pm \hv_0)\: \rd \bfz'+ \int_0^t \!\!\int_{\Upo_1}\!\! K_1 (\bfz,\bfz',t{-}s)r' \: (b_r \pm b_{\theta})\: \rd \bfz' \rd s. 
\end{equation*}
In view the rotational symmetry of $G_1$, the modified Green function is easily derived
\begin{equation*} 
G_2(r,\theta,r',\theta',t)=\frac{2}{\pi}\sum_{n,i=1}^{\infty}  \frac{n\sin(n(\theta-\theta')) }{(J_{n+1}(\kappa_{ni}))^2} \: J_{n}(\kappa_{ni} r) \: J_{n}(\kappa_{ni} r') \: \exp(-\kappa^2_{ni} \nu t).
\end{equation*}

By virtue of the asymptotic properties of the Bessel functions ($n\geq1$), we see that the kernel $K_2/r$ is well-behaved at the pipe centre, as
\begin{equation*}
n J_n(r) /r \sim {r^{n-1}}/{(2 (n-1)!)} \;\;\;\; \mbox{as}\;\;\; r \rightarrow 0. 
\end{equation*}
This evaluation ensures that the flow is regular at the axis of symmetry. It follows that both the integral kernels, $K_1$ and $K_2$, are bounded in the whole flow domain. By a standard algebraic manipulation to eliminate (say) $u_-$ from the equations (\ref{ie2}), the resulting equation becomes a linear Volterra equation of the second kind with a regular kernel. Its solutions can be established by successive approximations of resolvent kernels, see, for instance, Tricomi (1957). Without going into the technical details, we deduce that system (\ref{ie2}) is uniquely solvable, and the diffusion (\ref{apd}) admits solutions, $(\hu_r,\hu_{\theta}, \hu_z)$, for every given data (\ref{apic}) and buffer $\bfb$. 

As in the derivation of Poisson equation (\ref{ppoic}), system (\ref{apd}) always defines a solenoidal velocity field for any incompressible initial data, i.e., the last equation in (\ref{nsd}) holds at all times. Any momentum imbalance caused by the diffusions is restored, because the percolation gradients are fixed by the receptive non-linearities
\begin{equation} \label{apd2}
\partial_r \hat{q} = - \rho\alpha\beta f_r- \rho \beta b_r,\;\;\partial_{\theta} \hat{q} /r = - \rho\alpha\beta f_{\theta}-\rho \beta b_{\theta},\;\;\partial_z \hat{q} = - \rho\alpha\beta f_z - \rho \beta b_z,
\end{equation}
where the right-hand sides are known functions of $\bbu$ and $\hbu$.

If necessary, the pressure (\ref{ppcw}) can be calculated with extra effort. By virtue of the solutions of (\ref{apd2}), the Neumann condition is determined and denoted by
\begin{equation*} 
{\partial \hat{q}}/{\partial r}\big|_{r=1} = -\rho \beta b_r(1,\theta,z)=Q_0(\theta,z).
\end{equation*}
Both $\hq$ and $\bfb$ are assumed to decay at $|z|{\rightarrow} \infty$.
So the integral for the pressure is
\begin{equation} \label{pp}
\hat{q}(r,\theta,z) =\int_{\Upo_1} \!\!\!N_p(\bfz,\bfz') \:r'\:\hat{Q}_p(\bfz') \: \rd \bfz' + \int_{\real}\!\int_0^{2 \pi} 
\!\!\!N_p(\bfz,1,\theta',z') Q_0 \rd \theta' \rd z' + \hat{q}_0,
\end{equation}
where the Neumann kernel is governed by $ \Delta_0 N_p = - {\delta(r-r')\delta(\theta-\theta')\delta(z-z')}\big/r$, subject to homogeneous data. Starting with the expression for Dirac's delta,
\begin{equation*}
\delta(\theta-\theta')=\frac{1}{2 \pi}+ \frac{1}{\pi}\sum_{n=1}^{\infty}\cos(n(\theta-\theta'))=\frac{1}{2 \pi}\sum_{n=-\infty}^{\infty}\cos(n(\theta-\theta')),
\end{equation*}
it is straightforward to derive the product for the cross-section Neumann function
\begin{equation*}
\delta(r-r')\delta(\theta-\theta')/r=\frac{1}{\pi} \sum_{n=-\infty}^{\infty}\sum_{i=1}^{\infty}\frac{J_n(\sigma_{ni} r)J_n(\sigma_{ni} r')}{(1-n^2/\sigma^2_{ni}) (J_n(\sigma_{ni}))^2}\cos(n(\theta-\theta')),
\end{equation*}
where $\sigma_{ni}$ is the $i$th positive root of the function $J'_n(\sigma)=0$, so as to satisfy the boundary condition at the wall. By inspection, $N_p$ must have the form
\begin{equation*}
N_p(\bfz,\bfz')=\sum_{n=-\infty}^{\infty}\sum_{i=1}^{\infty}N_{ni}J_n(\sigma_{ni}r)\cos(n(\theta-\theta')).
\end{equation*}
In view of the pressure decay at large $|z|$, we determine $N_{ni}$ by solving the governing equation for $N_p$. The result is
\begin{equation*} 
 \begin{split}
	N_p= \frac{1}{\pi} \sum_{n=-\infty}^{\infty}\sum_{i=1}^{\infty}&\frac{J_n(\sigma_{ni} r)J_n(\sigma_{ni} r')}{(\sigma_{ni}-n^2/\sigma_{ni}) (J_n(\sigma_{ni}))^2} \; \times\\
	\quad & \hspace{20mm}\cos(n(\theta-\theta'))\exp\big(-\sigma_{ni}|z-z'|\: \big).
	\end{split}
\end{equation*}
In (\ref{pp}), the gauge $\hq_0$ does not affect the evaluation of gradients $\nabla\hq$.

There are no restrictions on the cycles of consecutive percolations. An interval, $[0,T]$, is sub-divided into $n$ equal parts. In each of these sub-intervals, we repeat our calculations for the local percolations. Over the last interval $[t_{n-1},t_n]$, the Hagen-Poiseuille profile $(\bbu,\baq)$ has companions
\begin{equation} \label{pipcn}
\bigg(\bbu+\sum_{k=1}^{n}\alpha^{(k)}\hbu^{(k)},\; \nabla \baq+\sum_{k=1}^{n}\nabla \hq^{(k)}/\beta^{(k)}\bigg)(x,y,z,t).
\end{equation}

The critical issue is that the complete Navier-Stokes equations admit non-linear as well as elementary solutions, though the class of the present solutions is not related to the most interesting phenomenon: turbulence. Our analysis demonstrates the fact that the companions decimate the parabolic Hagen-Poiseuille profile immediately after the commencement of the motion, as the percolations are multivalued at any time. Strictly, there does not exist one single identifiable characteristic velocity in the whole flow-field. Thus, the notion that turbulence in the pipe experiment of Reynolds (1883) is triggered by some form of finite-amplitude perturbation must be fundamentally flawed, as there are no well-defined Reynolds' numbers. Likewise, it is a truism to say that the laminar-turbulent transition in pipe flow spans a complex, non-linear, spatiotemporal process involving an enormous range of space and time scales, if such a rhetoric merely presumes the primitive particulars. In conclusion, flow instability and pleonastic bi-stability have no precise meanings.

\section{Couette motion between cylinders}

The flow between two rotating cylinders about a common axis has been studied extensively in fluid dynamics. In this section, we consider time-independent flows in the annulus between two rotating co-axial cylinders,
\begin{equation*}
\Upo_i:\;\;\;0< R_1 \leq r \leq R_2,\;\;\;0\leq \theta < 2 \pi,\;\;\; - \infty < z < + \infty,
\end{equation*}
where $R_1$ and $R_2$ are the radii of the inner and outer cylinders.
Denote their angular velocities by $\Omega_1$ and $\Omega_2$. The simplicity of the geometry suggests that the flow is characterised by the non-zero component $\bu_{\theta}=\bu_{\theta}(r)$ such that the incompressibility (\ref{divc}) is fulfilled. The governing equations (\ref{nsc}) reduce to
\begin{equation*}
\bu^2_{\theta}/r=\partial_r\bar{q}/\rho,\;\;\;\mbox{and}\;\;\;\nu (\partial_{rr} +\partial_r/r-1/r^2)\bu_{\theta} = 0,
\end{equation*}
respectively. The boundary conditions for real fluids ($\nu >0$) are
\begin{equation*}
\bu_{\theta}(R_1,\theta,z)=R_1\Omega_1,\;\;\;\bu_{\theta}(R_2,\theta,z)=R_2\Omega_2.
\end{equation*}
Algebra shows that the velocity is the sum of a rigid-body rotation and a viscous shear layer:
\begin{equation} \label{cf}
	\bbu=\big(\: 0,\;\;\; Ar+B/r,\;\;\;0\: \big),
\end{equation}
or its angular velocity:
\begin{equation} \label{cfa}
	\bar{\Omega}=\big(\: 0,\;\;\; A+B/r^2,\;\;\;0\: \big),
\end{equation}
where the integrating constants are expressed in symmetric forms
\begin{equation*}
A=(\Omega_2 R_2^2- \Omega_1 R_1^2)/(R_2^2-R_1^2),\;\;\;B=-(\Omega_2 - \Omega_1) R_2^2R_1^2/(R_2^2-R_1^2).
\end{equation*}
The mean pressure or its gradient is also known. Because of the no-slip condition,  the cylinders' walls drag the fluid into axi-symmetric rotational motion. 

\subsection{Fallacy of Rayleigh's instability criterion}

By formally putting $\nu=0$ in (\ref{nsc}), we get 
\begin{equation} \label{eulc}
	\frac{D_*}{D_*t}\tu_r-\frac{\tu^2_{\theta}}{r}= -\frac{1}{\rho}\frac{\partial \tp}{\partial r},\;\;\;\;\;
	\frac{D_*}{D_*t} (r\tu_{\theta}) = 0,\;\;\;\;\; \frac{D_*}{D_*t}\tu_z=-\frac{1}{\rho}\frac{\partial \tp}{\partial z}, 
\end{equation}
where $D_*/D_*t$ is the total derivative with axi-symmetry. The continuity reads,
\begin{equation} \label{idivc}
	{\partial \tu_r}/{\partial r}+{\tu_r}/{r}+{\partial \tu_z}/{\partial z}=0.
\end{equation}
The governing equations remain unchanged under $z$-translation. From the middle momentum of (\ref{eulc}), Rayleigh (1916) argued that the circulation, $r\tu_{\theta}$, or the angular momentum, $r^2\tilde{\Omega}$, is time-invariant. In the specific case of purely azimuthal flow, $\tu_{\theta} = \tu_{\theta}(r,t)$ where $\tu_r = \tu_z = 0$ or $\tu_r = 0$, $\tu_z = f(r,t)$, the $\theta$-momentum reduces to $\rd (r\tu_{\theta})/\rd t=0$. The kinetic energy per unit mass of the motion is $E_{\theta}=\tu^2_{\theta}/(2r^2)=(r^2\tilde{\Omega})^2/(2r^2)$. The stability of the motion depends on the sign of the discriminant,
\begin{equation} \label{ray}
\Upphi(r)=\frac{\rd} {\rd r} \big(r^2 E_{\theta}\big) = r^2{\tilde{\Omega}} \:\frac{\rd }{\rd r}\big(r^2 \tilde{\Omega}\big).
\end{equation}
By virtue of the conservation law, the constancy, $r^2\tilde{\Omega}(r,t)=r_0^2\tilde{\Omega}(r_0,0)$, does not necessarily imply stability or instability; we must restore to perturbation analyses (see below).

Rayleigh contemplated that the invariance of $r\tu_{\theta}$ still holds in three-dimensional motions with non-trivial radial and axial velocities, $\tu_r\neq0$ and $\tu_z\neq0$, provided the term $\tu^2_{\theta}/r$ in the first equation is balanced by a radial force. Then the whole motion would proceed as if $\tu_{\theta}$ were absent. His invariance postulation is imprecise without supplementary boundary conditions. Because of $z$ homogeneity, we consider only motions in planes normal to the axis, where $\tu_z$ is independent of $z$. From axi-symmetry, it is plain to see that continuity (\ref{idivc}) implies $\tu_r=\tilde{\alpha}(t)/r$ with any finite $\tilde{\alpha}$, and $\tilde{\alpha}(0)\neq0$. Then a solution of (\ref{eulc}) follows
\begin{equation} \label{iv0}
\tu_r=\frac{\tilde{\alpha}(t)}{r},\;\;\; r\tu_{\theta}=\exp\Big(-\frac{r^2}{2}+\int_0^t 
\tilde{\alpha}(s) \rd s\Big),\;\;\;\tu_z=g(r,t),
\end{equation}
as long as there are no boundary conditions. The pressure gradients are determined accordingly. Although the total derivative ${D_* (r\tu_{\theta})}/{D_*t}$ vanishes, it does not imply the time invariance of $r\tu_{\theta}$, because the circulation depends on the choice of function $\tilde{\alpha}(t)$ for fixed initial value $\tu_{\theta}(r,0)=\exp(-r_0^2/2)/r_0$. In an infinitely long hollow circular tube, $0\leq r \leq a$, solution (\ref{iv0}) is singular and inadmissible for flows of finite energy. Similarly, it does not exist wherever there is a boundary, on which $\tu_r$ vanishes. At least, the axi-symmetry and the condition $\tu_r|_{\partial \Upo}{=}0$ are incompatible propositions, unless $\tu_r=0$.

If the restriction on the axial symmetry is relaxed, Rayleigh's invariance on $r\tu_{\theta}$ is unavailable. Let us return to the inviscid flow in the hollow tube, with the boundary condition, $\tu(a,\theta,z,t)=0$. We find an arbitrary solution to the governing equations (\ref{nsc})-(\ref{divc}) ($\nu=0$),
\begin{equation} \label{iv1}
\tu_r = (r{-}a) \tilde{\alpha}(t)h(\theta),\;\;\tu_{\theta}=-(2r{-}a)\tilde{\alpha}(t) \!\int_0^{\theta}\!h(s)\rd s,\;\;\tu_z=Z(r,h(\theta))\tilde{\beta}(t).
\end{equation}
We assume $Z(r,h(\theta))$ are finite, or $Z=0$. Function $h$ is periodic, $h(\theta)=h(2 \pi m + \theta)$, where $m=0,\pm1,\pm2,\cdots$. The time functions, $\tilde{\alpha}(t)$ and $\tilde{\beta}(t)$, are arbitrary, smooth, and bounded, and $\tilde{\alpha}(0) \neq 0$, and $\tilde{\beta}(0) \neq 0$, to avoid trivial initial data. In fact, the inviscid momentum conservation is fully respected, once time-dependent $\tbu$ is specified, as each of $\nabla \tp$ in (\ref{eulc}) is fixed. 

For the corresponding exterior flow ($r\geq a$), a comparable solution is given by
\begin{equation} \label{ov1}
\tu_r = \frac{1}{r}(1-d_0)\tilde{\alpha}(t)h(\theta),\;\;\;\tu_{\theta}=-(k{+}1)(r-a)^k d_0\tilde{\alpha}(t) \!\int_0^{\theta}\!h(s)\rd s, 
\end{equation}
where the decay at large $r$ is controlled by $d_0=\exp\big(-(r{-}a)^{k+1}\big)$, $k \geq 1$. If we consider  the particular planar motion in real fluids, i.e., $\tu_z=0$ and $\nu>0$, the wall no-slip is satisfied by both $\tu_r$ and $\tu_{\theta}$. In other words, for given viscosity $\nu$, all the velocity terms in the planar Navier-Stokes equations (\ref{nsc}) are known so that the momentum conservation is maintained by the pressure gradients. Note that any suitable functions can be assigned to $\tilde{\alpha}(t)$ and $h(\theta)$ that, potentially, give rise to a family of disparate flow-fields. The novelty of the present example is that we do not have to rely on the use of the diffusive percolations to demonstrate the incompleteness of the primitive dynamics. 

\subsection*{Arbitrary perturbations}

In the annulus domain $\Upo_i$, the inviscid equations (\ref{eulc}) admit a steady mean flow $(\tbu,\tp)(r)$,
\begin{equation} \label{icf}
\tu_r=\tu_z=0,\;\;\;\tu_{\theta}=\tilde{V}(r)=r \tilde{\Omega}(r),\;\;\;\tp(r)=\int \big(\tilde{V}^2/r\big) \: \rd r,
\end{equation}
where the azimuthal velocity is an arbitrary function of $r$. Following the practice of perturbation theory, we consider the possibility that the mean flow is perturbed as a superposition
\begin{equation*}
\tu'_r, \; \tilde{V}(r)+\tu'_{\theta}, \; \tu'_z,\;\;\tp(r)+\tp'.
\end{equation*}
Now, the continuity is governed by
\begin{equation} \label{euldc}
{\partial \tu'_r}/{\partial r}+{\tu'_r}/{r}+({\partial \tu'_{\theta}}/{\partial \theta})/r + {\partial \tu'_z}/{\partial z}=0.
\end{equation}
The perturbation equations are written as
\begin{equation} \label{euld}
 {\partial \tp'}/{\partial r} =-\rho (\partial_t \tu'_r + f'_r),\;\;
 {\partial \tp'}/{\partial \theta} = -\rho r (\partial_t \tu'_{\theta} + f'_{\theta}),\;\;
 {\partial \tp'}/{\partial z}= -\rho (\partial_t \tu'_z + f'_z), 
 \end{equation}
where the full non-linearities are
\begin{align*}
 f'_r & =(\tbu'{\cdot}\nabla)\tu'_r + \tilde{\Omega}\:\partial \tu'_r/\partial \theta - {{\tu'}^2_{\theta}}/r - {2\tilde{\Omega} \tu'_{\theta}},\\ 
 f'_{\theta} & = (\tbu'{\cdot}\nabla)\tu'_{\theta} + \tilde{\Omega}\:\partial \tu'_{\theta}/\partial \theta  + \big(\rd \tilde{V}/\rd r+\tilde{\Omega}\big)\:\tu'_r + {\tu'_r \tu'_{\theta}}/{r}, \\
  f'_z & = (\tbu'{\cdot}\nabla)\tu'_z + \tilde{\Omega}\: \partial \tu'_z/\partial \theta.
\end{align*}
System (\ref{euldc})-(\ref{euld}) must be solved subject to the wall condition,
\begin{equation*}
\tu'_r(R_1,\theta,z,t)=\tu'_r(R_2,\theta,z,t)=0.
\end{equation*}
By observation, one set of the velocities is
\begin{equation} \label{iv2}
\begin{split}
\tu'_r &= (r-R_1)(R_2-r) \alpha(t) h(\theta),\\
\tu'_{\theta} & =\big(3r^2+2(R_2-R_1)r+R_1R_2\big)\alpha(t)\!\int_0^{\theta}\!h(s)\rd s,
\end{split}
\end{equation}
and $\tu'_z$ has the identical form of (\ref{iv1}). Another flow-field reads
\begin{equation} \label{iv3}
\tu'_r = \Big(\frac{1}{R_1}{-}\frac{1}{r}\Big) \Big(\frac{1}{r}{-}\frac{1}{R_2}\Big) \alpha(t) h(\theta),\;\;
\tu'_{\theta} =-\frac{1}{r}\Big(\frac{1}{r}{-}\frac{1}{R_1R_2}\Big)\alpha(t)\!\int_0^{\theta}\!h(s)\rd s.
\end{equation}
The sum of (\ref{iv2}) and (\ref{iv3}) constitutes additional solutions. Given an inviscid Couette mean and the perturbations, equations (\ref{euld}) indicate that the gradients are exactly known. Note that these solutions also satisfy the usual linearised equations, where the non-linear terms in the functions, $f'_r,f'_{\theta}$, and $f'_z$, are neglected. It becomes clear that the gradient $\nabla \tp'$ needs to be updated to reflect the approximations due to linearisation. The point is that our finite-energy velocity fields can be superimposed on any claimed solutions of (\ref{euld}) which establish unstable properties of the mean. It follows that there is no instability as postulated by criterion (\ref{ray}). 

\subsection{Circular shear between infinitely long cylinders}

Couette motion (\ref{cf}) or its inviscid counterpart (\ref{icf}) of the unbounded cylindrical gap is translation-invariant in the $z$-direction. Hence the wave-trains of the specific form (\ref{dst3}) are not supported in these flows. 

To establish percolations (\ref{ap}), we solve equations (\ref{apd}) subject to initial data similar to those in (\ref{apic}), and the non-slip condition on $\partial \Upo_i$.
Instead of Green's function $G_0$ in (\ref{g0}), we need 
\begin{equation*} 
 \begin{split}
	g_0&= \frac{\pi}{4} \sum_{i=1}^{\infty} F_0(\lambda_{ni}) E_0(\lambda_{ni} r) E_0(\lambda_{ni} r')\exp(-\lambda^2_{ni} \nu t) \; + \\
	\quad & \hspace{15mm} \frac{\pi}{2} \!\sum_{n,i=1}^{\infty} \!\cos(n(\theta-\theta'))F_n(\lambda_{ni}) E_n(\lambda_{ni} r) E_n(\lambda_{ni} r') \exp(-\lambda^2_{ni} \nu t),
	\end{split}
\end{equation*}
where $J_n$ and $Y_n$ stand for the Bessel functions of order $n$, and $F_n$ denotes the normalisation factor,
\begin{equation} \label{norf}
	F_n(\lambda_{ni})=\frac{\lambda_{ni}^2(J_n(\lambda_{ni} R_2))^2}{(J_n(\lambda_{ni} R_1))^2-(J_n(\lambda_{ni} R_2))^2},
\end{equation}
and $E_n$ the expansion functions
\begin{equation*} 
	E_n(\lambda_{ni} r)=J_n(\lambda_{ni} R_1)\;Y_n(\lambda_{ni} r) - Y_n(\lambda_{ni} R_1)\;J_n(\lambda_{ni} r).
\end{equation*}
In every summation over $n=0,1,2,\cdots$, $\lambda_{ni}$ is the $i$th positive root of the equation
\begin{equation} \label{g1s}
	J_n(\lambda R_1)\;Y_n(\lambda R_2) - Y_n(\lambda R_1)\;J_n(\lambda R_2)=0.
\end{equation}
Similarly, we obtain the analogous system to (\ref{ie2}) in terms of a regular kernel $K_2$, where $G_2$ is replaced by 
\begin{equation*} 
	g_2=\frac{\pi}{2} \!\sum_{n,i=1}^{\infty} \!n \sin(n(\theta-\theta')) \; F_n(\lambda_{ni}) \; E_n(\lambda_{ni} r) \; E_n(\lambda_{ni} r') \; \exp(-\lambda^2_{ni} \nu t).
	\end{equation*}
Furthermore, the integral kernel to replace $G_1$ is
\begin{equation*} 
	g_1=\frac{\pi}{2} \!\sum_{n,i=1}^{\infty} \!\cos(n(\theta-\theta')) \; F_n(\lambda_{ni}) \; E_n(\lambda_{ni} r) \; E_n(\lambda_{ni} r') \; \exp(-\lambda^2_{ni} \nu t).
\end{equation*}
We obtain the gradients directly from the law of momentum conservation
\begin{equation*} 
\begin{split}
\frac{\partial \hq}{\partial r} &= - \rho\alpha\beta \:  \Big( \alpha(\hbu{\cdot}\nabla)\hu_r-\alpha\frac{\hu^2_{\theta}}{r}+\Big(A{+}\frac{B}{r^2}\:\Big)\Big(\: \frac{\partial \hu_r}{\partial \theta} - 2 \hu_{\theta}\: \Big)\Big) - \rho \beta b_r, \\
\frac{\partial \hq}{\partial \theta} &= - \rho\alpha\beta \:  \Big( \alpha \: r (\hbu{\cdot}\nabla)\hu_{\theta} + \Big(\:Ar{+}\frac{B}{r}\:\Big) \Big(\:\frac{\hu_{\theta}}{\partial \theta} +\hu_r\:\Big) + \Big(\:A{-}\frac{B}{r^2}\:\Big)\hu_r\Big)- \rho \beta b_{\theta}, \\
\frac{\partial \hq}{\partial z} & = - \rho\alpha\beta \: \Big(\alpha(\hbu{\cdot}\nabla)\hu_z +\Big(\:A{+}\frac{B}{r^2}\:\Big)\frac{\partial \hu_z}{\partial \theta} \Big)- \rho \beta b_z. 
\end{split}
\end{equation*}
In short, the companions have the form (\ref{pipcn}).

\subsection{Perturbations to Taylor's linearised solution}

If the disturbances are assumed to take a symmetric form about the axial direction, we drop all terms $\partial/\partial \theta$ in the equations of motion. Taylor (1923) considered the superimposed flow of the Couette profile (\ref{cf}) and linear axi-symmetric disturbances, $\bbu+\hbu^T, \bp+\hq^T$. The perturbation flow is governed by 
\begin{equation} \label{lcf}
\begin{split}
\partial_t \hu^T_r & - 2 \big(A+B/r^2\big) \hu^T_{\theta}  = \nu \Delta_1^* \hu^T_r-\frac{1}{\rho}\frac{\partial \hq^T}{\partial r}, \\
\partial_t \hu^T_{\theta} + 2A \hu^T_r & = \nu \Delta_1^* \hu^T_{\theta},\;\;\;
\partial_t \hu^T_z = \nu \Delta^*_0 \hu^T_z- \frac{1}{\rho}\frac{\partial \hq^T}{\partial z},
\end{split}
\end{equation}
and by the continuity
\begin{equation} \label{lca}
\nabla^*{\cdot}\hbu^T=\frac{\partial \hu^T_r}{\partial r}+\frac{\hu^T_r}{r}+\frac{\partial \hu^T_z}{\partial z}=0.
\end{equation}
The Laplacians $\Delta_1^*=\Delta_0^*-1/r^2$, and 
$\Delta_0^* = {\partial^2}/{\partial r^2} + ({\partial}/{\partial r})/r + {\partial^2}/{\partial z^2}$.
By making use the expressions in (\ref{ppoic}) for axi-symmetric flow, we find
\begin{equation} \label{pac}
\Delta_0^* \hq^T = 2 \rho \: \Big(\: \big( \: A+{B}/{r^2} \: \big)\:{\partial \hu^T_{\theta}}/{\partial r} + \big(\:A-B/r^2\:\big)\: {\hu^T_{\theta}}/{r}\:\Big).
\end{equation}
In the narrow-gap approximation, Taylor sought the solutions of (\ref{lcf}) and (\ref{lca}) by expanding the disturbances in terms of the eigen-functions,
\begin{equation} \label{eigfn}
	\hu^T_r=\tu_r (r) \cos(k z),\;\;\; \hu^T_{\theta}=\tu_{\theta}(r) \cos(k z),\;\;\; \hu^T_{z}=\tu_z(r) \sin(k z),
\end{equation}
with boundary conditions $\hbu^T=0$ at $r=R_1,R_2$. On the end-plates, $z=0,2 \pi/k$, the inviscid condition applies to $\hu^T_z=0$; $\hu^T_r$ and $\hu^T_{\theta}$ are periodic in the axial direction.

Let the annulus be unit height,
$\Upo_a:R_1 \leq r \leq R_2, 0\leq \theta < 2 \pi, 0 \leq z \leq 1$.
Denote the corresponding Taylor's steady solutions by
\begin{equation} \label{tflow}
\bbu+\hbu^T=\bfU=(U,\;\;V,\;\;W)(r,z),\;\;\;\bp+\hq^T=\bar{P}(r,z).
\end{equation}
Now we consider a set of perturbations 
\begin{equation} \label{pc2}
(\bfU+\alpha \bfv,\;\;\bar{P}+\tp/\beta)(r,z,t).
\end{equation}
On substituting (\ref{pc2}) into governing equations (\ref{nsc}), retaining the axi-symmetry and linearising, we get the following dynamic equations:
\begin{equation} \label{lcf1}
\begin{split}
\partial_t v_r +(\bfU{\cdot}\nabla^s)v_r + (\bfv{\cdot}\nabla^s)U - 2 {V v_{\theta}}/{r} & = \nu \Delta_1^* v_r-({\partial \tp}/{\partial r})/(\rho\alpha\beta), \\
\partial_t v_{\theta} +(\bfU{\cdot}\nabla^s)v_{\theta} + (\bfv{\cdot}\nabla^s)V + {U v_{\theta}}/{r}+{V v_r}/{r} & = \nu \Delta_1^* v_{\theta},\\
\partial_t v_z +(\bfU{\cdot}\nabla^s)v_z + (\bfv{\cdot}\nabla^s)W& = \nu \Delta^*_0 v_z-({\partial \tp}/{\partial z})/(\rho\alpha\beta),
\end{split}
\end{equation}
where the gradient operator $\nabla^s={\partial}/{\partial r}+{\partial}/{\partial z}$, and $\nabla^*{\cdot}\bfv=0$. The no-slip condition on $\partial \Upo_a$ reads
\begin{equation} \label{abc}
\bfv(R_1,z,t)=\bfv(R_2,z,t)=0;\;\;\;\bfv(r,0,t)=\bfv(r,1,t)=0.
\end{equation}
We are interested in the initial-boundary value problem (\ref{lcf1}) subject to data
\begin{equation*} 
\bfU_s=(U_s,\;\;V_s,\;\;W_s)(r,z,t=0).
\end{equation*}

For any axi-symmetric buffer $\bfb=(b_r,b_z)(\bfz,t)$, the percolations are defined by
\begin{equation*} 
\partial_t v_r-\nu \Delta_1^* v_r=b_r,\;\;\;\mbox{and}\;\;\;\;\partial_t v_z - \nu \Delta^*_0 v_z=b_z.
\end{equation*}
The solution for the $z$-component is readily found to be
\begin{equation} \label{vz2}
	v_z(\bfz,t)=\int_{\Upo_a}\!\! K_a(\bfz,\bfz',t) r' W_s(\bfz') \: \rd \bfz' + \int_0^t \!\! \int_{\Upo_a}\!\! K_a(\bfz,\bfz',t{-}s) r' b_z(\bfz',s) \: \rd \bfz' \rd s,
\end{equation}
where $K_a(\bfz,\bfz',t)=H_f(z,z',t)M_0(r,r',t)$, and $M_0$ is given by 
\begin{equation*} 
	\frac{\pi}{4} \sum_{i=1}^{\infty} \; F_0(\lambda_{0i}) \; E_0(\lambda_{0i} r) \; E_0(\lambda_{0i} r') \; \exp(-\lambda^2_{0i} \nu t), 
\end{equation*}
(see (\ref{norf}) to (\ref{g1s}) at $n=0$). Similarly,
\begin{equation} \label{vr2}
	v_r(\bfz,t)=\int_{\Upo_a}\!\! H_a(\bfz,\bfz',t)r' U_s(\bfz') \: \rd \bfz' +\int_0^t \!\! \int_{\Upo_a}\!\! H_a(\bfz,\bfz',t{-}s) r' b_r(\bfz',s) \: \rd \bfz' \rd s,
\end{equation}
where $H_a(\bfz,\bfz',t)=H_f(z,z',t)M_1(r,r',t)$, and 
\begin{equation*} 
	M_1= \frac{\pi}{4} \sum_{i=1}^{\infty} \; F_1(\lambda_{1i}) \; E_1(\lambda_{1i} r) \; E_1(\lambda_{1i} r')\; \exp(-\lambda^2_{1i} \nu t).
\end{equation*}
In comparison with (\ref{pac}), the forcing term on the pressure, $-\Delta_0^* \tp /(\rho\alpha\beta)$, has a more complicated expression
\begin{equation*}
\Big(\frac{\partial}{\partial r}+\frac{1}{r}\Big)\Big( (\bfU.\nabla^s)v_r + (\bfv.\nabla^s)U - 2 {V v_{\theta}}/{r} \Big) +  \frac{\partial}{\partial z}\Big( (\bfU.\nabla^s)v_z + (\bfv.\nabla^s)W\Big).
\end{equation*}
We give this equation for the sake of completeness. The Neumann boundary values are inhomogeneous because the derivatives of $\bfb$ do not vanish on the boundary unless the condition, $\bfb|_{\bdy_a}=0$, is imposed. In practice, we need only the gradients. To balance the diffusive $v_r$ and $v_z$, we equate the $z$-gradient to 
\begin{equation*} 
- \rho \alpha \beta \:\big(\: U \partial_r v_z + W \partial_z v_z + v_r \partial_r W + v_z \partial_z W  \:\big) - \rho \beta b_z.
\end{equation*}
From (\ref{lcf}), the kernel of operator $\partial_t{-}\nu \Delta^*_1$ is simply $M_1$. Then the integral for $v_{\theta}$ is
\begin{equation} \label{iet}
\begin{split}
v_{\theta}(\bfz,t)&= \bigg( \int_{\Upo_a} \!\! M_1 V_s \: \rd \bfz'+ \int_0^t \!\!\int_{\Upo_a} \!\! M_1 \big(\: (\bfv.\nabla')V + V v_r/r \:\big) \rd \bfz' \rd t' \bigg)\\
& \hspace{5mm}+ \int_0^t \!\!\int_{\Upo_a} \!\! (M_1 U/r) \: v_{\theta} \:\rd \bfz' \rd t' + \int_0^t \!\!\int_{\Upo_a} \!\! M_1 (U \partial_r  + W \partial_z) \:v_{\theta} \:\rd \bfz' \rd t' \\
& = F_a(\bfz,t) + \int_0^t \!\!\int_{\Upo_a} \!\! M(\bfz,\bfz',t{-}s) \: v_{\theta}(\bfz',s) \:\rd \bfz' \rd s,
\end{split}
\end{equation}
where the revised kernel, 
$M = M_1 U / r  - {\partial (M_1U)}/{\partial r} - {\partial (M_1W)}/{\partial z}$, in view of the bounded $U$ and $W$ on $\bdy_a$, and the no-slip $v_{\theta}$. It follows that $v_{\theta}$ is known, as equation (\ref{iet}) is a linear Volterra integral equation of the second kind with a bounded regular kernel. Since $\bfv$ has been found, the radial gradient satisfies
\begin{equation*}
{\partial \tp}/{\partial r} = - \rho \alpha \beta \: \big(\: U \partial_r v_r + W \partial_z v_r + v_r \partial_r U + v_z \partial_z U -2V v_{\theta}/r \:\big)- \rho \beta b_r,
\end{equation*}
for the compatibility of the first equation in (\ref{lcf1}). Instead of 
the linearisation and the axi-symmetry, the analogous full non-linear system exists, and its solution can be sought in the same way as we did for the case of infinitely long cylinders. The multiplicity of percolations (\ref{pc2}) implies that the perturbative proposition for the annulus flow between co-axial rotating cylinders, as stipulated by the linearised boundary-value problem (\ref{lcf})-(\ref{lca}), loses its potential for prediction. 

\section{Remark on flow similarity}

If we write equation (\ref{ns}) in dimensionless form
\begin{equation*} 
	\partial_t \bfu^0 + (\bfu^0 {\cdot} \nabla) \bfu^0  = - \nabla p^0/\rho + \Delta \bfu^0/Re,\;\;\; \nabla{\cdot}\bfu^0=0, 
\end{equation*}
for $Re=UL/\nu$, we can still exercise the flow partition as follows:
\begin{equation*}
\partial_t \bfu^0 -Re^{-1}\Delta \bfu^0= \bfb^0;\;\;\;\nabla p^0=-\rho (\bfu^0 {\cdot} \nabla) \bfu^0- \rho \bfb^0.
\end{equation*}
Since the partition is driven by arbitrary dimension-free buffer $\bfb^0$, the $Re$-specified primitive formulation does not admit unique solutions. As there can be no consensus on well-defined velocity fields, the idea of flow similarity is self-contradictory.

The analysis of the Hagen-Poiseuille flow shows that there is not a single characteristic velocity anywhere inside the pipe. The only sensible choice for a velocity scale lies in its initial value. Yet, different motions can be generated in a pipe apparatus (diameter $d$). In practice, there are no particular reasons that the following profiles cannot be initiated by dedicated pressure-differentials:
\begin{enumerate}
\item A parabola, $u_z(r,t{=}0)=U_0(1-r^2)$;
\item An oscillation, $u_z(r,t{=}0)=U_0\big(\sin(5\pi(1/2-r))+1\big)/2$;
\item A precipice, $u_z(r,t{=}0)=U_0\tanh(\pi^2(1-r))$;
\item A slide, $u_z(r,t{=}0)=U_0\big(1-\exp(-\pi^2(1-r)^2)\big)$,
\end{enumerate}
assuming $\nabla{\cdot}\upomega_0=0$. It is logical to select the maximum velocity $U_0$ as a characteristic scale. It is plausible that the local flows over downstream distances $\sim O(5d {-} 10d)$ must be distinct in these four scenarios. The implication is that similarity descriptions of these motions by one dimensionless number, $U_0d/\nu$, are too simplistic to accommodate the diversity of the initial data. Thus, the nature of flow similarities departs from its original formalism. Possibly, a logical way forward is to investigate classifications of starting flows in order for non-dimensional groups to be consequential.

\subsection*{Vorticity setting}

In principle, fluid dynamics is governed by the conservation law of angular momentum. Then the vorticity, $\upomega=\nabla {\times} \bfu=(\xi,\eta,\zeta)$, is of paramount importance (Helmholtz 1858). Dynamically, the vorticity evolves according to  
\begin{equation} \label{vort}
	\partial_t \upomega - \nu \Delta \upomega = (\upomega {\cdot} \nabla) \bfu  - (\bfu {\cdot} \nabla )\upomega,
\end{equation}
which is to be solved subject to data (\ref{vtic}) and (\ref{vtbc}).
The compatibility, $\bfu = (\Delta)^{-1} (\nabla{\times}\upomega)$, implies that $\bfu$ is driven by $\upomega$. The fact is that only the vorticity is {\it a priori} bounded, and it has an invariant property under Galilean transform,
$\upomega=\nabla{\times}\bfu=\nabla{\times}(\bfu-{\bf c})$. Note that a partition of the equation ceases to work, as the pressure is truly dependent, see (\ref{ppoi}). Briefly, the continuum description of flow evolution, including turbulence, is tantamount to the interactions of {\it both} non-linear terms. They are responsible for instigating multitudinous vortices. 

If we take out the common dimension factor $[\texttt{s}^{-2}]$, equation (\ref{vort}) becomes
\begin{equation} \label{nvt} 
	\partial_t \upomega^0 + (\bfu^0 {\cdot} \nabla )\upomega^0 -(\upomega^0 {\cdot} \nabla) \bfu^0 =  \nu^0 \Delta \upomega^0,
\end{equation}
where the viscosity $\nu^0$ is relative to unity. This normalised equation shows that the denominator $[\texttt{s}^2]$ is a measure of enstrophy. Hence, a non-dimensional parameter may be found by reference the kinematic viscosity, $\mu/\rho$, to a representative circulation $\omega L^2$.
Plausibly, dynamic similarities may be understood for individual categories of initial data in similar geometries. 

The rate of dissipation of kinetic energy per unit mass (denoted by $\varepsilon_E$) has the dimension $[\texttt{(m/s)}^2/\texttt{s}]$. As implied in (\ref{vort}), the instantaneous dissipation occurs across the whole range of the vorticity eddies, which measure the local strains and shears. On the contrary, the aggregate effects of all the vortices assemble the local velocity. In the absence of a representative velocity, this rate of dissipation may be understood as the energy induced by the local eddies, $\varepsilon_E \sim O\big(\vec{\bfu}{\cdot}\vec{\upomega}\big)$, where $\vec{\bfu}=(u^2,v^2,w^2)$, and $\vec{\upomega}=(\sqrt{\xi^2},\sqrt{\eta^2},\sqrt{\zeta^2})$ are the energy and eddy vectors respectively. This estimate is just one of the possible results by dimensional analysis. An analogy to the familiar $O\big(|\bfu|^3/L_{Turb}\big)$ is an approximate, $O\big(L^2_{Turb} |\upomega|^3\big)$, which is interpreted as the circulation of local enstrophy. Given vorticity $\upomega$, the governing equations (\ref{ns}) assert that the energy dissipates according to $\varepsilon_E= - \nu |\upomega|^2$,
which defines a {\it glocal} dynamic process, depending on the viscosity at all times. 

\section{Helmholtz's decomposition}

The reason why the formulation of the primitive variables is ill-advised is due to the fact that the role of the vorticity has not been accounted for. To analyse the vorticity dynamics, we must first consider the well-known vector decomposition. In his paper on vorticity theory, Helmholtz (1858) shows that any smooth velocity field in the whole Euclidean space $\rr$ can be decomposed into the gradient of a scalar potential and the curl of a vector potential:
\begin{equation} \label{hm}
\bfu = \nabla \phi + \nabla {\times} \uppsi, 
\end{equation}
consult \S147-\S148 of Lamb (1932); \S20 of Sommerfeld (1950); \S 1.5 and \S2.3 of Morse \& Feshbach (1953). Since $\bfu$ is incompressible, we have $\Delta \phi=0$. Taking curl of the decomposition, vector calculus gives $\Delta \uppsi = \nabla(\nabla{\cdot}\uppsi)-\nabla{\times}\bfu$. The usual practice is to set $\nabla{\cdot}\uppsi=0$. 
In smooth flows decaying at infinity, Laplace's and Poisson's equations are well-posed, and hence the decomposition is necessarily unique. 

In application, there are two scenarios that call for special attention: (a) Generalised flow fields; (b) Incomplete decomposition. First, the velocity has a potential,
\begin{equation} \label{hm2}
	\bfu = \nabla \phi, \;\;\;\;\; (\nabla{\times}\uppsi=0)
\end{equation}
where the vanishing of the curl is assumed to hold. This is a matter of steady solutions for the Euler's equations. Second, in the vorticity dynamics, it is paramount to recover the velocity field from the vorticity,
\begin{equation} \label{hm4}
	\upomega = \nabla{\times}\bfu,
\end{equation}
where, once again, it is assumed that a gradient ($\nabla \phi$) does not contribute. We calculate $\bfu$ by solving either 
\begin{equation} \label{uv}
	\Delta \bfu(\bfx) = - \nabla{\times}\upomega(\bfx),\;\;\; \bfx \in \rr.
\end{equation}
or $\Delta\uppsi=-\upomega$, where $\bfu=\nabla{\times}\uppsi$. In the latter approach, the gradient is missing in view of the full decomposition (\ref{hm}). The key is that (\ref{hm2}) and (\ref{hm4}) are {\it incomplete} decompositions as some properties of the vector fields are neglected. There is a compatibility issue in these formulas. 

In flow-fields having distributional derivatives, the continuity is now interpreted as
\begin{equation} \label{dw2}
	\int_{\rr}(\nabla{\cdot}\bfu)\;\tilde{v}\;\rd \bfx = - \int_{\rr}\bfu \cdot \nabla\tilde{v} \;\rd \bfx = 0
\end{equation}
for every test function $\tilde{v} \in C_0^{\infty}(\rr)$. The right integral makes sense for $\bfu \in L^1_{loc}$. Likewise, we interpret irrotational vectors as generalised derivatives as follows:
\begin{equation} \label{dw4}
	\int_{\rr}(\nabla{\times}\bfu)\cdot\tilde{\bfv}\;\rd \bfx = - \int_{\rr}\bfu \cdot(\nabla{\times}\tilde{\bfv}) \;\rd \bfx = 0
\end{equation}
for every vector test function $\tilde\bfv$. There must exist irrotational test functions. An example is constructed in
\begin{equation} \label{dtf}
\tilde{\bfv}_i(\bfx)=\left\{
\begin{aligned}
&\;\big( E_0,\;\; E_0, \;\; E_0 \big), &\;\;\; s < 1, \\
&\;\;0, &\;\;\; s \geq 1, \\
\end{aligned}
\right.
\end{equation}
where $s=\sqrt{x^2+y^2+z^2+2xy+2yz+2zx}$, and
\begin{equation} \label{e0}
E_0(\bfx)=C\exp\bigg( - \frac{1}{1 - s^2}\bigg),
\end{equation}
for normalisation constant $C$.

Let us start from the full Helmholtz decomposition for vorticity,
\begin{equation} \label{vhm}
	\upomega=\nabla{\times}\bfu + \nabla \varphi, \;\;\;\;\;(\nabla{\cdot}\bfu=0).
\end{equation}
Multiplying (\ref{vhm}) by $\nabla \tilde{v}$ for test function $\tilde{v} \in C_0^{\infty}$ and integrating by parts, we arrive at
\begin{equation}
	\int_{\rr} \tilde{v} \; \big(\nabla{\cdot} \upomega - \Delta \varphi \big) \;\rd \bfx = 0.
\end{equation}
Since this integral vanishes for every test function $\tilde{v}$, we must have,
\begin{equation} \label{dv2}
	\Delta \varphi = \nabla{\cdot} \upomega, \;\;\; \mbox{a.e.}
\end{equation}
For unique (vector) stream function $\uppsi$ in (\ref{hm}), we put $\uppsi=\uppsi'+\nabla\phi'+{\bf c}$, where $\Delta \phi'=-\nabla{\cdot}(\uppsi'+{\bf c})$, and $\bf c$ is any irrotational constant. Then $\bfu$ can be recovered by solving
\begin{equation} \label{vpo}
\Delta \uppsi = - \upomega,\;\;\; \nabla{\cdot}\uppsi=0,
\end{equation}
assuming $\uppsi$ decays. Conversely, the vector potential is given by 
\begin{equation} \label{vp}
\uppsi(\bfx) = \frac{1}{4 \pi}\int_{\rr} \frac{1}{|\bfx-\bfy|}\;\upomega(\bfy)\; \rd \bfy=\frac{1}{4 \pi}\int_{\rr} \frac{1}{\;|\bfy|\;}\;\upomega(\bfx{+}\bfy)\; \rd \bfy.
\end{equation}
The Biot-Savart law shows that gauge choice, $\nabla{\cdot}\uppsi=0$, holds everywhere, provided vorticity is solenoidal, $\nabla{\cdot}\upomega=0$. By Gauss's divergence theorem (cf. (\ref{dv})), $\upomega$ decays at infinity. Hence potential $\varphi$ must be harmonic:
\begin{equation} \label{poi}
	\Delta \varphi =  0, \;\;\; \varphi \rightarrow 0, \;\;\; |\bfx| \rightarrow \infty,
\end{equation}
where the decay is necessary to avoid harmonic functions which are unbounded at infinity, such as, $\varphi = C_1 x^2 +C_2 y^2 - (C_1+C_2)z^2$. The decomposition components are orthogonal as
\begin{equation} \label{vvn}
\begin{aligned}
	\int_{\rr} \nabla{\times}\bfu \cdot \nabla \varphi \; \rd \bfx & = \big[ \varphi \nabla{\times}\bfu \big]_{[\bfx]}- \int_{\rr} \varphi \;(\nabla {\cdot} \nabla{\times}\bfu) \; \rd \bfx \\
	& = \big[ \bfu\cdot \nabla\varphi \big]_{[\bfx]} - \int_{\rr} \bfu \cdot (\nabla{\times} \nabla \varphi) \; \rd \bfx = 0,
\end{aligned}
\end{equation}
where the last two integrands vanish identically by vector identities.

Taking dot product of (\ref{vhm}) with vector test function $\nabla{\times}\tilde{\bfv}$, and integrating, we get the integral condition
\begin{equation}
	\int_{\rr} \tilde{\bfv} \cdot \big(\Delta \bfu + \nabla{\times}\upomega \big) \;\rd \bfx = 0,
\end{equation}
which shows that velocity $\bfu$ is recovered by solving elliptic equation (\ref{uv}).

Conversely, given a solenoidal vector, $\nabla{\cdot}\upomega=0$, where $\upomega$ decays at large $|\bfx|$, how do we decompose this vector field? 
Let $\bff$ be a non-vanishing vector on $\rr$ which differs from $\upomega$. Furthermore, we impose the condition $\bff \rightarrow 0$ at large distances. Vector analysis shows that
\begin{equation} \label{v2}
\Delta \bff = \nabla(\nabla{\cdot}\bff) - \nabla{\times} \nabla {\times} \bff.
\end{equation}
By curl of test function $\tilde{\bfv}$, this differential equation is formulated as 
\begin{equation} 
\begin{aligned}
\int_{\rr} \nabla{\times}\tilde{\bfv} \cdot \Delta \bff \; \rd \bfx &= - \int_{\rr} \tilde{\bfv} \cdot\Delta (\nabla{\times}\bff)\; \rd \bfx \\
 & = -\int_{\rr} \tilde{\bfv}\cdot \nabla{\times} \nabla(\nabla{\cdot}\bff)\; \rd \bfx + \int_{\rr} \tilde{\bfv} \cdot(\nabla{\times}\nabla{\times}\nabla{\times}\bff)\; \rd \bfx \\ & = - \int_{\rr} \tilde{\bfv} \cdot \Delta \bfg\; \rd \bfx,
\end{aligned}
\end{equation}
where we introduce function $\bfg$, and $\nabla{\cdot}\bfg=0$, and $\nabla{\times}\nabla{\times}\bfg=-\Delta \bfg$. Since it is solenoidal, vector $\bfg$ is better characterised as $\bfg|_{|\bfx| \rightarrow \infty} \rightarrow 0$. This integral condition reduces to
\begin{equation} \label{v4}
\int_{\rr} \tilde{\bfv} \cdot \Delta (\nabla{\times}\bff -\bfg)\; \rd \bfx = 0,
\end{equation}
that shows that $\bfg$ differs from $\nabla{\times}\bff$ within multiples of a harmonic function. We denote it by $\bfh$, where $\nabla{\cdot}\bfh=0$, $\nabla{\times}\bfh=0$, giving $\bfg=\nabla{\times}\bff+\bfh$ a.e. As in potential flows, $\bfh$ is understood as a generalised function according to (\ref{dw2}) and (\ref{dw4}). Without loss of generality, we choose a scalar potential $\varphi$ such that $\nabla \varphi = \bfh$, and hence $\varphi$ satisfies (\ref{poi}). 

In identity (\ref{v2}), we replace $\nabla{\times}\bff$ by $\bfg$. The solution can be expressed by
\begin{equation}
\bff(\bfx)= \frac{1}{4 \pi}\int_{\rr}\frac{1}{|\bfx-\bfy|}\nabla(\nabla{\cdot}\bff)\;\rd \bfy + \frac{1}{4 \pi}\int_{\rr} \nabla_{\bfy}\Big(\frac{1}{|\bfx-\bfy|}\Big){\times}\bfg(\bfy) \; \rd \bfy.
\end{equation}
For every test function $\tilde{v} \in C_0^{\infty}(\rr)$, we evaluate
\begin{equation}
\int_{\rr} \bff\cdot \nabla\tilde{v} \; \rd \bfx = - \int_{\rr} \tilde{v}\;\nabla{\cdot}\bff \; \rd \bfx =  (I_1 + I_2)/(4 \pi),
\end{equation}
where
\begin{equation}
\begin{aligned}
I_1 & = - \int_{\rr} \tilde{v}\; \int_{\rr} \nabla_{\bfx}{\cdot}\Big(\frac{1}{|\bfx-\bfy|}\Big)\nabla(\nabla{\cdot}\bff)\;\rd \bfy \; \rd \bfx\\
& = \int_{\rr} \tilde{v}\; \bigg( \int_{\rr} \Big( \frac{\bfx-\bfy}{|\bfx-\bfy|^3}\Big)\cdot\nabla(\nabla{\cdot}\bff)\;\rd \bfy \bigg) \; \rd \bfx,
\end{aligned}
\end{equation}
and
\begin{equation}
I_2=-\int_{\rr} \tilde{v}\; \bigg( \int_{\real^3 }\nabla_{\bfx}{\times}\nabla_{\bfy}\Big(\frac{1}{|\bfx-\bfy|}\Big)\cdot \bfg(\bfy) \; \rd \bfy \bigg) \; \rd \bfx =0,
\end{equation}
where the last equality is established due to the $\bfx{-}\bfy$ symmetry in the Newtonian potential.
A quick recap tells us that the divergence is determined by a homogeneous integro-differential equation with an integrable kernel: 
\begin{equation} \label{v6}
g(\bfx) = -\frac{1}{4 \pi}\int_{\rr} {\bf K}(\bfx-\bfy)\cdot \nabla g(\bfy)\;\rd \bfy,
\end{equation}
for shorthand $g(\bfx)=\nabla{\cdot}\bff(\bfx)$. Denote the Fourier transform of $g$ by
\begin{equation} \label{ft}
	{\mathscr F}(g)(\bfk)={\hat g}(\bfk)=\frac{1}{(\sqrt{2 \pi})^3}\int_{\rr}\exp( - \ri \;\bfk \cdot \bfx) \;g(\bfx) \;\rd \bfx.\;\;\;(\ri = \sqrt{-1})
\end{equation}
We transform equation (\ref{v6}), as the kernel is a convolution type. Re-arranging the result leads to algebraic formula
\begin{equation} \label{w2}
{\hat g}(\bfk)\;\big(\;1 + \ri \;\bfk \cdot{\mathscr F}({\bf K})\sqrt{{\pi}/{2}} \;\big)=0,
\end{equation}
which is valid for arbitrary wave-number vector $\bfk$. The only consistency is ${\hat g} \equiv 0 $. It follows that its inverse
\begin{equation} \label{w4}
\nabla{\cdot}\bff(\bfx)= \frac{1}{(\sqrt{2 \pi})^3} \int_{\rr}\exp(\ri \;\bfk \cdot \bfx) \;{\hat g}(\bfk) \;\rd \bfk = 0,
\end{equation}
identically. 

In summary, the vector field $\bff$ characterised by differential constraints, (\ref{v2}) and (\ref{v4}), must be incompressible and tends to zero at infinity. The essential reason for its existence lies in the solenoidal property of $\bfg$. It follows that vector $\bfg$ can be identified with vorticity $\upomega$, and $\bff$ with velocity $\bfu$. In fact, we have recovered the Helmholtz decomposition: $\upomega=\nabla{\times}\bfu+\nabla\varphi$, where the kinematics is understood in the generalised sense. In particular, the presence of harmonic $\varphi$ is a direct consequence of $\nabla{\cdot}\upomega=0$. Vector algebra suggests that the velocity may be non-unique, as $\nabla{\times}\bfu=\nabla{\times}(\bfu+\nabla{q}+{\bf c})$, where $q$ is any scalar function. This arbitrariness is insignificant, as what is being decomposed is the vorticity vector, i.e., the velocity is computed from (\ref{uv}) for given vorticity. 

In fact, our discussion of full Helmholtz's decomposition links to the properties of Laplace's equation, as defined in (\ref{poi}). Multiplying it by test function $\tilde{v} \in C_0^{\infty}(\real^3)$ and integrating, we get
\begin{equation*} 
	\int_{\rr}\Delta \varphi \;\tilde{v} \; \rd \bfx = - \int_{\rr} \nabla \varphi \cdot \nabla \tilde{v} \; \rd \bfx = \int_{\rr}\varphi \; \Delta \tilde{v} \; \rd \bfx = 0,
\end{equation*}
where the weak function $\varphi$ is smooth on $\rr$ in view of Weyl's lemma on the interior regularity (Weyl, 1940). By (\ref{dw2}) and (\ref{dw4}), the theorem reiterates the physics that there are no sources generating vorticity in the `interior' of the whole space, hence ruling-out discontinuous derivatives. Nevertheless, Weyl's proof does not include the `boundary': the origin and infinity. 
 
For $\varphi \in C^2(\rr)$, Green's identity gives 
\begin{equation*} 
	\int_{\rr} \varphi \Delta \varphi \; \rd \bfx = - \int_{\rr} \big(\nabla \varphi \cdot \nabla\varphi \big) \; \rd \bfx = - \int_{\rr} |\nabla \varphi|^2 \; \rd \bfx =0.
\end{equation*}
Thus, the only smooth solution is trivial $\varphi \equiv 0$. This fact has long been established, see p.490 of Helmholtz (1858) for regular simply-connected domains, and \S35 of Poincar\'e (1893) for $\rr$. The textbook concept, $\upomega\equiv\nabla{\times}\bfu$, must have been made on the presumption that flow-fields of $\upomega,\bfu$ are the classical solutions of their respective governing equations. 

If the kinematic solutions are considered in the functional spaces other than $C^{\lambda}(\rr)$, $\lambda \geq 2$, well-documented harmonic functions show that either they are unbounded at infinity (an example is $\varphi = Ax + Byz + Cxyz$), or they possess a singularity at the origin, as a source or a sink or a doublet. The former case is rejected on the ground of finite energy. It follow that Laplace's equation holds on $\rr{-}{\bfx_1}$ only. Specifically, the gradient is well-behaved at locations away from $\bfx_1$,
\begin{equation*}
\nabla\varphi = - \frac{\kappa}{4 \pi} \frac{(\bfx-\bfx_1)}{\;|\bfx-\bfx_1|^3},\;\;\;\;\;\bfx\neq\bfx_1,
\end{equation*}
for arbitrary $\kappa$. Because of the singularity, Helmholtz's decomposition on the vorticity (or the velocity) breaks down at $\bfx=\bfx_1$. If $\upomega$ becomes singular at location $\bfx=\bfx_s$, this singularity has no significance, as $\nabla \varphi$ can be translated to $\bfx_s=\bfx_1$.

In view of (\ref{vpo}) and (\ref{vp}), the gauge choice, $\nabla{\cdot}\uppsi=0$, implies $\nabla{\cdot}\upomega=0$. It is evident that vector identity, $\nabla{\times}\nabla{\times}\uppsi=-\upomega$, is inherently solenoidal by formula $\nabla{\cdot}\nabla{\times}{\bf A}=0$. From the sufficiency point of view, an arbitrary gauge of $\uppsi$ cannot be verified by the Biot-Savart integral in the absence of $\upomega$ regularity. Let us fix a constant or any continuous function $f(t)=C=\nabla{\cdot}\uppsi\neq0$. Solution (\ref{vp}) is now self-contradictory. 

\section{Vorticity bounds}

By vector identity, $(\bfg{\cdot}\nabla)\bfg = \nabla{\times}\bfg{\times}\bfg  + \nabla(|\bfg|^2)/2$, the momentum equation is rewritten as
\begin{equation} \label{ns2}
\partial_t \bfu  - \nu \Delta \bfu - \bfu \times (\nabla {\times} \bfu) = - \nabla \chi,
\end{equation}
where $\chi = \chi(\bfu,p)=p/\rho + |\bfu|^2/2$, is known as the total pressure.

The principle of conservation of angular momentum holds in fluid motions. Applying this law to control volume analysis, vorticity dynamics can be derived. Analytically, the pressure can be eliminated from system (\ref{ns2}). Vector algebra shows that compatibility of solenoidal $\bfh(\bfx,t)$ requires $\Delta \bfh = - \nabla{\times}\nabla{\times}\bfh$ at every given instant. Consider test functions, $\tilde{\bfv}(\bfx,t) \in C_0^{\infty}(\rr){\times}(0,\infty)$. We take dot product on momentum (\ref{ns}) by $\nabla{\times}\tilde{\bfv}$, integrate the result over space and time to obtain
\begin{equation} \label{vt2}
\begin{split}
\int_{\real}\int_{\rr}\bfu \cdot \partial_t (\nabla{\times}\tilde{\bfv}) \; \rd \bfx \rd t + \nu \int_{\real}\int_{\rr} &(\nabla{\times}\bfu) \cdot (\nabla{\times}\nabla{\times}\tilde{\bfv}) \; \rd \bfx \rd t \\
& + \int_{\real}\int_{\rr} \bfu{\times}(\nabla {\times} \bfu) \cdot (\nabla{\times}\tilde{\bfv})\; \rd \bfx\rd t = 0,
\end{split}
\end{equation}
where $p$ is eliminated from integration by parts as 
\begin{equation} 
\int_{\real}\int_{\rr} p \; \nabla{\cdot}(\nabla{\times}\tilde{\bfv}) \; \rd \bfx\rd t = 0.
\end{equation}
The continuity also holds in the weak sense
\begin{equation} \label{ct2}
  \int_{\real}\int_{\rr}(\nabla {\cdot} \bfu) \;\nabla{\times}\tilde{\bfv}\; \rd \bfx\rd t=-\int_{\real}\int_{\rr} \bfu \;\nabla{\cdot}(\nabla{\times}\tilde{\bfv})\; \rd \bfx\rd t = 0.
\end{equation}
We call the quantity, $\nabla {\times} \bfu$, vorticity, denoted by $\upomega$, the components of which are $(\xi,\eta,\zeta)$. Vorticity describes the local rotation of fluid elements.
Equations, (\ref{vt2})-(\ref{ct2}), make sense for $\bfu \in L^1 \cap L^2$, $\upomega \in L^1 \cap L^2$. Now suppose that $\bfu$ is sufficiently regular, so that we can perform further integrations. We have
\begin{equation} 
\int_{\real}\int_{\rr}\tilde{\bfv}\cdot \big(\;\partial_t (\nabla{\times}\bfu) - \nu \Delta (\nabla{\times}\bfu) + \nabla{\times}(\nabla {\times} \bfu){\times}\bfu\;\big)\; \rd \bfx\rd t = 0,
\end{equation}
for every test function $\tilde{\bfv}$. Thus the vorticity dynamics is written in (\ref{vort}).
In inviscid flows, the contribution from the viscous term is ignored by setting $\nu\equiv0$. The advantage of the vorticity dynamics is evident as the pressure gradient now exerts no mechanical moment or couple on fluid elements. Furthermore, the vorticity field inherits the incompressibility constraint $\nabla{\cdot}\upomega=0$.  

\subsection*{Preliminaries}

Given a multi-index of order $n \;(\geq 1)$, $\alpha=(\alpha_1,\; \alpha_2,\; \cdots, \alpha_n)$ with positive $\alpha$'s, denote $|\alpha|=\alpha_1+\alpha_2+\cdots +\alpha_n$. Denote the $\alpha$ derivative by
\begin{equation*}
D^{\alpha}f(\bfx) = \frac{\partial^{|\alpha|}}{\partial x_1^{\alpha_1} \partial x_2^{\alpha_2} \cdots \partial x_n^{\alpha_n}} f(\bfx)= \partial_{x_1}^{\alpha_1}\partial_{x_2}^{\alpha_2} \cdots \partial_{x_n}^{\alpha_n} f(\bfx).
\end{equation*}
For two functions $f(\bfx)$ and $g(\bfx)$, derivatives of their products can be found by Leibniz's rule,
\begin{equation} \label{lb}
D^{\alpha}(fg)=\sum_{\alpha \geq \beta} \;\frac{\alpha!}{\beta!\:(\alpha-\beta)!}\; (D^{\beta}f)\; (D^{\alpha-\beta}g),
\end{equation}
where $\alpha$ and $\beta$ are two multi-indices.

To allow for possible unbounded vorticity or velocity, we will discuss our solutions and their derivatives in the weak or generalised sense. For integers, $1 \leq p \leq \infty$ and $k \geq 0$, we denote the Sobolev space by
\begin{equation} \label{sb}
W^{k,p}(\rr) = \big\{ v \in L^p(\rr): \;\;\;D^{\alpha} v \in L^p(\rr), \;\;\;\mbox{for all} \;\; \alpha \in {\mathbb N},\; |\alpha| \leq k \big\},
\end{equation}
where $D^{\alpha}v$ are the distributional derivatives, see, for instance, Chapter 3 of Adams \& Fournier (2003) for well-known properties. Specifically,  
\begin{equation} \label{ww}
W_0^{k,p}(\rr) = W^{k,p}(\rr).
\end{equation}
The Sobolev norm, $\|\cdot\|_{W^{k,p}(\rr)}$, is defined by
\begin{equation}
\begin{aligned}
\|v\|_{W^{0,p}(\rr)} & = \|v\|_{L^p(\rr)}=\Big( \int_{\rr} | \:v\: |^p\; \rd \bfx\Big)^{1/p},\;\;\;\;\;(1\leq p < \infty)\\
\|v\|_{W^{k,p}(\rr)} & = \bigg( \sum_{|\alpha|\leq k} \big\|D^{\alpha}v \big\|^p_{W^{0,p}(\rr)}\bigg)^{1/p}, \\
\|v\|_{W^{0,\infty}(\rr)} & = \|v\|_{L^{\infty}(\rr)} =  \underset{\bfx \in \rr}{\mbox{ess\,sup}} \; |v(\bfx)|,\\
\|v\|_{W^{k,\infty}(\rr)} & = \max_{|\alpha| \leq k }\|D^{\alpha}v\|_{L^{\infty}(\rr)}.
\end{aligned}
\end{equation}

All bold-faced letters are vector qualities. On occasion, we drop the space dependence and write $\|\bfu\|_{L^p(\rr)}$ as $\|\bfu\|_{L^p}$. Frequently, we have to deal with integrals with integrand $f \nabla g$. We introduce notation for the boundary term,
\begin{equation*}
\int_{\rr} \!\! f \partial_x g\; \rd \bfx = - \int_{\rr} \!\! g \partial_x f \; \rd \bfx+\lim_{x \rightarrow \pm \infty} \int_{\real^2} f g\; \rd y \rd z =-\int_{\rr} \!\! g \partial_x f \; \rd \bfx + \big[ fg \big]_{[x]}.
\end{equation*}
The meanings of $[fg]_{[y]}$ and $[fg]_{[z]}$ are clear by cyclic permutation. All three can be summed and extended in $[\bfg{\cdot}\bff]_{[\bfx]}$ as well as in $[\bfg{\times}\bff]_{[\bfx]}$.

\subsection*{The continuity}

The incompressibility, $\nabla{\cdot}\bfu_0=\nabla{\cdot}\bfu=0$, is a finite sum of three reals. It follows that every $\partial u_i/\partial x_i$ must be finite. There is a natural number $N_1$ such that every one of these three derivatives can be bounded
\begin{equation}
\bigg|\frac{\partial u_i}{\partial x_i}\bigg| < N_1 \; \sum_{n=1}^3\;\; \sup_{\bfx \in \rr} \bigg|\frac{\partial u_n(\bfx,t{=}0)}{\partial x_n}\bigg|.
\end{equation}
The conservation of mass actually imposes strong restrictions on the magnitude of the velocity or the vorticity. Conceptually, it is difficult to comprehend how the continuity would be maintained somehow when $\bfu$ becomes unbounded, as fluid material would have been annihilated\textemdash an absurd physics. 

At $t=0$, the postulation of bounded data (\ref{vtic}) enables us to fix a useful bound
\begin{equation}
A^*_0=N_2\bigg(\;\sup_{\bfx \in \rr}\big\|D^{\beta}\xi_0\big\|_{L^1} + \sup_{\bfx \in \rr}\big\|D^{\beta}\eta_0\big\|_{L^1} + \sup_{\bfx \in \rr}\big\|D^{\beta}\zeta_0\big\|_{L^1}\;\bigg),
\end{equation}
for $\beta=0,1,2,\cdots$, natural number $N_2$.

\subsection{Properties of non-linearity}

Denote the total velocity and vorticity by 
\begin{equation}
U(\bfx,t) = u(\bfx,t) {+} v(\bfx,t) {+} w(\bfx,t),\;\;\;\mbox{and},\;\;\;
\Omega(\bfx,t) = \xi(\bfx,t) {+} \eta(\bfx,t) {+} \zeta(\bfx,t),
\end{equation}
respectively. The inviscid vorticity equation ($\nu=0$ in (\ref{vort})) is
\begin{equation} \label{k2}
	\partial_t \upomega = (\upomega \cdot \nabla) \bfu  - (\bfu \cdot \nabla )\upomega.
\end{equation}
Consider the x-component $\xi$. Integrating over space, we obtain
\begin{equation} \label{k4}
\begin{aligned}
\frac{\rd}{\rd t}\int_{\rr} \xi \; \rd \bfx = \big[ u &\xi\big]_{[x]} + \big[ u \eta\big]_{[y]} + \big[ u \zeta\big]_{[z]}- \int_{\rr} (\nabla{\cdot}\upomega)\; \bfu \; \rd \bfx \\
& - \big[ u \xi \big]_{[x]} - \big[ v \xi \big]_{[y]} -\big[ w \xi \big]_{[z]} + \int_{\rr} (\nabla{\cdot}\bfu)\; \upomega \; \rd \bfx=0,  
\end{aligned}
\end{equation}
where all the boundary terms vanish by virtue of the decay (\ref{vtbc}). We repeat these steps for the other two components, $\eta$ and $\zeta$. This result states that the total torque on fluid elements equals to zero. The dynamics equation for the total vorticity reads
\begin{equation} \label{k6}
	\partial_t \Omega = (\upomega \cdot \nabla) U  - (\bfu \cdot \nabla )\Omega.
\end{equation}
Use of the total vorticity is advantageous in scenarios where a pair of counter-rotating vortices involves highly concentrated strong shears. Integration yields
\begin{equation} \label{k8}
\begin{aligned}
\frac{\rd}{\rd t}\int_{\rr} \Omega \; \rd \bfx = \big[ U &\xi\big]_{[x]} + \big[ U \eta\big]_{[y]} + \big[ U \zeta\big]_{[z]}- \int_{\rr} (\nabla{\cdot}\upomega)\; U \; \rd \bfx \\
& - \big[ u \Omega \big]_{[x]} - \big[ v \Omega \big]_{[y]} -\big[ w \Omega \big]_{[z]} + \int_{\rr} (\nabla{\cdot}\bfu)\; \Omega \; \rd \bfx=0.  
\end{aligned}
\end{equation}
Every component, as well as the total, satisfy an invariance principle
\begin{equation}
\int_{\rr}\upomega(\bfx,t)\; \rd \bfx =\int_{\rr}\upomega_0(\bfx)\;\rd\bfx,\;\;\;\;\;\int_{\rr}\Omega(\bfx,t)\;\rd\bfx = \int_{\rr}\Omega_0(\bfx)\;\rd\bfx.
\end{equation}
But 
\begin{equation*}
\int_{\rr} \xi \;\rd \bfx= \int_{\rr} \max(\xi,0) \;\rd \bfx + \int_{\rr} \min(\xi,0) \;\rd \bfx,
\end{equation*}
where the last sum is finite and hence each of the integrals must be bounded. Thus, $\xi$ is integrable, as are $\eta$ and $\zeta$. We come to the summability of the field
\begin{equation} \label{m4}
\xi, \; \eta, \; \zeta \in L^1(\rr),\;\;\;\mbox{and},\;\;\;\Omega \in L^1(\rr).
\end{equation}
By the Archimedean property for two positive reals,
\begin{equation} \label{m6}
\|\upomega\|_{L^1(\rr)} < A^*_0.
\end{equation}
As $\upomega$ is a member of the equivalent class $W^{0,1}$, it can be well approximated by equivalent class of $C_0^{\infty}$ functions on $\rr$, in view of (\ref{ww}). In the sequence of approximations, each dies out like $ \sim |\bfx|^{-\varepsilon} \;(\varepsilon>0)$ or $ |\bfx|^{-2}$ as a polynomial in $\nabla{\cdot}\upomega=0$ near infinity. In the absence of precise decay rates, which depend on the initial data, the best we can infer is
\begin{equation} \label{m8}
\upomega \rightarrow 0, \;\;\;\;\; |\bfx| \rightarrow \infty.
\end{equation}
We have
\begin{equation} \label{j2}
\upomega \in L^{\infty}(\Upo_T)\times W^{(0,1)}_0(\rr).
\end{equation}

\subsection*{Weak derivatives}

Denote test function $\tilde{v}(\bfx,t) \in C_0^{\infty}(\rr){\times}(0,\infty)$
Multiplying the $x$-component of vorticity equation (\ref{vort}) by $D^{\alpha} \tilde{v}$, $\alpha=1,2,3,\cdots$ and integrating, we obtain 
\begin{equation} \label{vtx}
	\int_{\real}\int_{\rr} \big( \; \xi \; (\partial_t + \nu \Delta ) D^{\alpha}\tilde{v} - (\upomega \;u - \bfu \; \xi) \cdot \nabla D^{\alpha}\tilde{v} \;  \; \big)\; \rd \bfx \rd t =0,
\end{equation}
in view of $\nabla{\cdot}\bfu=0$ and $\nabla{\cdot}\upomega=0$. Together with the other components, we see that equations (\ref{vtx}) make sense for $\bfu, \upomega \in L^1 \cap L^2$. For the inviscid case $\alpha=1$, we carry out integration by parts, and acquire a conservation law for the variation of the torque
\begin{equation}
	\int_{\real}\int_{\rr} \tilde{v} \;\big(\; \partial_t (D \upomega) - D (\; (\upomega \cdot \nabla) \bfu  - (\bfu \cdot \nabla )\upomega )\; \big)\; \rd \bfx=0,
\end{equation}
for every test function. Hence, the differential equation for vorticity derivative reads
\begin{equation} \label{dvt}
	\partial_t (D \upomega) = D \big(\; (\upomega \cdot \nabla) \bfu  - (\bfu \cdot \nabla )\upomega \; \big).
\end{equation}
Consider $D=\partial_x$ first. Integration over space shows that
\begin{equation}
\frac{\rd}{\rd t}\int_{\rr} \partial_x \upomega \; \rd \bfx  = \big[ (\upomega \cdot \nabla) \bfu  - (\bfu \cdot \nabla )\upomega \big]_{[x]}=0,  
\end{equation}
where use is made of estimate (\ref{m8}). In particular, no use has been made of the assumed regularity on the right-hand function in (\ref{dvt}). Analogous operations on the total vorticity admit the invariance
\begin{equation}
\frac{\rd}{\rd t}\int_{\rr} \partial_x \Omega \; \rd \bfx  = \big[ (\upomega \cdot \nabla) U  - (\bfu \cdot \nabla )\Omega \big]_{[x]}=0.  
\end{equation}
Repeating these integrations with $D=\partial_y$ and $D=\partial_z$ in turn, we assert that
\begin{equation} \label{dvt2}
\xi, \; \eta, \; \zeta \in W_0^{1,1}(\rr),\;\;\;\mbox{and}\;\;\;\Omega \in W_0^{1,1}(\rr),
\end{equation}
respectively. By the rationale leading to (\ref{m8}), the rate of change in vorticity must tend to zero at infinity,
\begin{equation} \label{om}
D\upomega \rightarrow 0, \;\;\;\;\; |\bfx| \rightarrow \infty.
\end{equation}
The upgrade of (\ref{j2}) is obvious.

Simple algebra gives $D\upomega=\nabla{\times}D\bfu$ and $\nabla{\cdot}D\bfu=D(\nabla{\cdot}\bfu)=0$, implying the velocity-vorticity compatibility, $\Delta D\bfu=-\nabla{\times}D\upomega$. The velocity inversion is done by
\begin{equation} \label{bs2}
D\bfu(\bfx) = -\frac{1}{4 \pi}\int_{\rr} \Big(\frac{\bfx-\bfy}{|\bfx-\bfy|^3}\Big) \times D\upomega(\bfy)\; \rd \bfy,
\end{equation}
where the boundary term vanishes in view of results (\ref{dvt2}) and (\ref{om}). We denote
\begin{equation*}
\mbox{supp} \;D\upomega \subset \overline{B_1(\bfx_0,R_1)}, \;\;\;\;\;R_1=R_1(\bfu_0) \gg 1.
\end{equation*}
For $|\bfx| \gg |\bfy| \sim O(R_1)$, the denominator of the kernel can be approximated as
\begin{equation}
|\bfx-\bfy| = \sqrt{\bfx{\cdot}\bfx}\;\sqrt{ 1 {-} 2\; (\bfx{\cdot}\bfy)/(\bfx{\cdot}\bfx) {+} (\bfy{\cdot}\bfy)/(\bfx{\cdot}\bfx)} \approx |\bfx|.
\end{equation}
As $|\bfx|\rightarrow \infty$,
\begin{equation} \label{bc2}
D\bfu(\bfx) \sim \frac{\bfx}{\;|\bfx|^3\;} {\times} \! \int_{B_1} \!\! D\upomega(\bfy)\; \rd \bfy \sim  \frac{\bfx}{\;|\bfx|^3\;} {\times} \!\int_{B_0} \!\! D\upomega_0(\bfy)\; \rd \bfy \rightarrow 0,
\end{equation}
because the last integral is bounded. The support of initial data $\upomega_0$ is denoted by $B_0$. This decay implies $\bfu \equiv 0$ at infinity. In addition, $\upomega|_{|\bfx| \rightarrow \infty} \equiv 0$ by (\ref{om}). This kinematic state is a consequence of dynamics (\ref{k2}), the boundary condition (\ref{vtbc}), and the solenoidal properties. Specifically, no claim has been made on the regularity of $D\bfu$ in (\ref{bs2}). Applying this procedure of estimation to $\Delta \bfu=- \nabla{\times}\upomega$ as a consistency check, we recover the $\bfu$ boundary condition for integrable vorticity. The physics dictates the existence of the stagnation, because the vicinity at infinity must remain undisturbed, and cannot be accessed by any localised incompressible flow. One implication is that all higher derivatives of the velocity have to vanish at large far distances. However, we shall proceed with our inductive reasoning. 

The weak differentiation on (\ref{dvt}) gives us an equation for the second-order derivatives or $\alpha=2$ in (\ref{vtx}),
\begin{equation} 
	\partial_t D^2 \upomega = D \big( D \big(\; (\upomega \cdot \nabla) \bfu  - (\bfu \cdot \nabla )\upomega \; \big)\big).
\end{equation}
For $\partial_t \partial_x D \upomega$, explicitly, the boundary term from the right-hand side,
\begin{equation*}
\big[ D ((\upomega \cdot \nabla) \bfu)  - D ((\bfu \cdot \nabla )\upomega) \big]_{[x]} =0,
\end{equation*}
by deductions (\ref{om}) and (\ref{bc2}). Note that we have no knowledge of the `interior regularity'. In summary, $D^2\upomega \in W_0^{0,1}$, and $D^2\upomega$ decays. 

For $\alpha \geq 2$, since $\nabla{\cdot}D^{\alpha}\bfu=D^{\alpha}(\nabla{\cdot}\bfu)=0$, it follows that
\begin{equation} 
\Delta D^{\alpha}\bfu = - \nabla \times D^{\alpha} \upomega.
\end{equation}
Inverting the Laplacian yields
\begin{equation} 
D^{\alpha}\bfu(\bfx)=-\frac{1}{4 \pi} \int_{\rr} \Big(\frac{\bfx-\bfy}{|\bfx-\bfy|^3}\Big) \times D^{\alpha}\upomega(\bfy)\; \rd \bfy,
\end{equation}
provided that $D^{\alpha}\upomega$ is integrable and tends to zero,
\begin{equation} \label{omk}
 D^{\alpha}\upomega \rightarrow 0, \;\;\;\;\; |\bfx| \rightarrow \infty.
\end{equation}
Following the analysis leading to (\ref{bc2}), we establish
\begin{equation} \label{bc4}
D^{\alpha}\bfu \rightarrow 0, \;\;\;\;\; |\bfx| \rightarrow \infty,
\end{equation}
which are valid regardless of the differentiability of the velocity at locations $\bfx \ll \infty$.
We take derivatives of order $\alpha$ on (\ref{dvt}) to get
\begin{equation} 
	\partial_t (D^{\alpha+1} \upomega) = D \big\{ D^{\alpha} \big(\; (\upomega \cdot \nabla) \bfu  - (\bfu \cdot \nabla )\upomega \; \big)\big\},
\end{equation}
as we wish to show that $D^{\alpha+1}\upomega$ is integrable. In fact, we only need to consider the term in the outer brackets at infinity. Each term in $D^{\alpha}(\cdot)$ can be expanded according to Leibniz's rule (\ref{lb}). In each of the two expansions, every resulting term consists of products of derivatives of the vorticity and velocity. All these derivatives are equal to zero as $x,y,z \rightarrow \pm \infty$, in view of (\ref{omk}), (\ref{bc4}), as well as the decays of the lower-orders. Induction shows
\begin{equation} \label{n7}
\frac{\rd}{\rd t}\int_{\rr} D^{\alpha} (\xi,\; \eta,\; \zeta)\; \rd \bfx =0,\;\;\;\;\;\alpha=1,2,3,\cdots.
\end{equation}
An alternative interpretation of this outcome suggests that the integrations and differentiations commute in view of the attenuation at large distances, even though the smoothness of the flow in finite spheres is undetermined. Consequently, the bound analogous to (\ref{m6}) is given by
\begin{equation}
\|D^{\alpha}\upomega\|_{L^1(\rr)} < A^*_0,
\end{equation}
and hence
\begin{equation} \label{d6}
\upomega \in W_0^{k,1}(\rr),\;\;\;\;\;k=1,2,3,\cdots.
\end{equation}
By virtue of the Sobolev embedding theorem (see, for instance, Theorem 4.12 of Adams \& Fournier, 2003), we arrive at a space-wise bound 
\begin{equation} \label{m10}
\upomega \in C^2_B(\rr),
\end{equation}
which, {\it a posteriori}, justifies the assumption made in (\ref{nsc}), as the partial Helmholtz decomposition $\upomega = \nabla{\times}\bfu$
enshrines the definition of vorticity at every given instant. The point is that it is best to recover the velocity from the vorticity by solving (\ref{uv}), 
$\Delta \bfu (\bfx) = - \nabla{\times} \upomega (\bfx)$. Remarkably, the integrability of vorticity ($L^1$-norms) and the `boundary' decay (\ref{vtbc}) are the essential ingredients. The space-time bound (\ref{j2}) is improved to
\begin{equation} \label{j4}
\upomega \in L^{\infty}(\Upo_T)\times C^2_B(\rr).
\end{equation}
Additional information of the velocity is known according to the well-developed theory of elliptic partial differential equations (see, for instance, Gilbarg \& Trudinger, 1998). For the time being, it is sufficient to notice
\begin{equation} \label{m12}
\bfu \in C^2_B(\rr).
\end{equation}

For H\"older continuous $\upomega \in C^{0,\lambda}(\rr)$, $0 \leq \lambda <1$, it is well-established that $\nabla\bfu \in C^1(\rr)$, which is regular enough to justify the omission of the harmonic component in $\upomega=\nabla{\times}\bfu+\nabla\varphi$. The implication is that our integrations of the vorticity field in the Sobolev space $W^{k,1}$ are very effective but are not necessarily optimum. As $C^0$ is embedded in $W^{3,1}$, the next regularity upgrade results in $\upomega \in W^{4,1}$ which misses all fractional values of $\lambda$.
 
Taking divergence on Euler's momentum,
\begin{equation}
\partial_t \bfu = - (\bfu{\cdot}\nabla)\bfu - \nabla p/\rho,
\end{equation}
gives us Poisson's equation for the pressure
\begin{equation} 
\Delta p(\bfx) = - \rho \; \nabla \cdot \big((\bfu{\cdot}\nabla)\bfu \big) = - \rho \; \nabla \cdot \big( (\bfu{\cdot}\nabla)\bfu \big) = -\rho N_0(\bfu),
\end{equation}
where every term on the right-hand side has been shown to be continuous and bounded, see (\ref{m10}) and (\ref{m12}). Consequently, the pressure gradient is bounded in space. In particular, the gradient tends to zero at large distances, giving rise to the gauge pressure at infinity. 

\subsection{Vorticity dynamics for real fluids}

Let the vorticity equation,
\begin{equation} \label{mvt}
	\partial_t \upomega - \nu \Delta \upomega = (\upomega \cdot \nabla) \bfu  - (\bfu \cdot \nabla )\upomega \equiv {\bf V}(\bfx,t).
\end{equation}
This dynamics is to be solved subject to data (\ref{vtic}) and decay (\ref{vtbc}). The total vorticity satisfies the scalar equation
\begin{equation} \label{ntv}
	\partial_t \Omega - \nu \Delta \Omega = (\upomega \cdot \nabla) U  - (\bfu \cdot \nabla )\Omega.
\end{equation}
By Duhamel's principle, equation (\ref{mvt}) can be transformed into an equivalent integral equation
\begin{equation} \label{ie3}
\upomega(\bfx,t)=\int_{\rr}  H(\bfx{-}\bfy,t) \upomega_0(\bfy)\; \rd \bfy + \int_0^t \!\! \int_{\rr}  H(\bfx{-}\bfy,t{-}s) \;{\bf V}(\bfy,s) \; \rd \bfy \rd s,
\end{equation}
where $H$ is the diffusion kernel 
\begin{equation} \label{dk}
H(\bfx,t)= \frac{1}{\;(4 \pi t)^{3/2}}\; \exp \bigg( - \frac{\;|\bfx|^2}{4 t} \bigg).
\end{equation}
Integration of (\ref{ie3}) over space defines a constraint
\begin{equation}
\int_{\rr} \upomega(\bfx,t) \; \rd \bfx = \int_{\rr} \upomega_0(\bfx)\; \rd \bfx + \int_0^t \! \Big( \int_{\rr} \! {\bf V}(\bfx,s) \; \rd \bfx \Big) \rd s.
\end{equation}
The last integral equals to zero, as analysed in the preceding paragraphs. Thus the invariance,
\begin{equation}
\frac{\rd}{\rd t} \int_{\rr} \upomega \; \rd \bfx = 0,\;\;\;\;\;(\nu > 0)
\end{equation}
holds, whereby viscous $\upomega$ satisfies bound (\ref{m4}). Evidently, the above derivation applies equally well to the total vorticity governed by (\ref{ntv}). 

In view of (\ref{vtx}), a sequence of differential equations also works for weak $\upomega$-derivatives,
\begin{equation}
	(\partial_t  - \nu \Delta )D^{\alpha }\upomega = D^{\alpha}{\bf V}(\bfx,t).
\end{equation}
For every $\alpha$, the analogous equation to (\ref{ie3}) reads
\begin{equation}
D^{\alpha}\upomega(\bfx,t)=\int_{\rr} \! H(\bfx{-}\bfy,t) D^{\alpha}\upomega_0(\bfy)\; \rd \bfy + \int_0^t \!\! \int_{\rr} \! H(\bfx{-}\bfy,t{-}s) \;D^{\alpha} {\bf V}(\bfy,s) \; \rd \bfy \rd s,
\end{equation}
and hence
\begin{equation} 
\int_{\rr} D^{\alpha}\upomega(\bfx,t) \; \rd \bfx = \int_{\rr} D^{\alpha}\upomega_0(\bfx)\; \rd \bfx + \int_0^t \! \Big( \int_{\rr} D^{\alpha}{\bf V}(\bfx,s) \; \rd \bfx \Big) \rd s.
\end{equation}
Without repeating the analysis on the integral of $D^{\alpha}{\bf V}$, we assert that bounds, such as (\ref{n7}) and (\ref{d6}), are also valid for the viscous equation. 

In summary, we derive {\it a priori} bound 
\begin{equation} \label{d8}
|\upomega|,\; |\nabla \bfu| < M_0(\bfu_0)< \infty, \;\;\;\;\; (\bfx \in \rr).
\end{equation}
These space-wise estimates are effective at every given instant. 

\subsection{Enstrophy and energy}

The rate of change of enstrophy can be computed as follows. Taking dot product of equation (\ref{vort}) by $\upomega$ and integrating the result, we see that
\begin{equation*}
\frac{1}{2}\frac{\rd}{\rd t}\int_{\rr}|\upomega|^2 \; \rd \bfx + \nu \int_{\rr}|\nabla\upomega|^2 \; \rd \bfx = \int_{\rr}\upomega\cdot (\upomega{\cdot}\nabla)\bfu \; \rd \bfx - \int_{\rr}\upomega\cdot (\bfu{\cdot}\nabla)\upomega \; \rd \bfx,
\end{equation*}
where we have made use of the vorticity decay on the second term on the left. In view of the solenoidal velocity and vorticity, the last term on the right vanishes as a result of integration by parts. The remaining term is estimated as
\begin{equation*}
\int_{\rr}|\upomega\cdot (\upomega{\cdot}\nabla)\bfu| \; \rd \bfx \leq \|\upomega\|_{L^{\infty}} \|\upomega \nabla\bfu\|_{L^1} \leq \|\upomega\|_{L^{\infty}} \|\upomega\|^2_{L^2},
\end{equation*}
by applying H\"older's inequality. Thus, the evolution is governed by
\begin{equation} \label{es}
\frac{\rd}{\rd t}\|\upomega\|^2_{L^2} < 2 M_0 \|\upomega\|^2_{L^2},
\end{equation}
the solution of which bounds the enstrophy,
\begin{equation}
\|\upomega(\cdot,t)\|^2_{L^2} < \|\upomega_0\|^2_{L^2} \exp\big( 2 M_0 t\big).
\end{equation}
The exponential is likely an over-estimate. Importantly, the law of energy conservation holds at all times $t>0$,
\begin{equation} 
\frac{1}{2}\int_{\rr} \! \rho|\bfu(\bfx,t)|^2 \; \rd \bfx + \mu \int_0^t \!\! \int_{\rr} \! |\upomega(\bfx,t)|^2\; \rd \bfx = \frac{1}{2}\int_{\rr} \! \rho|\bfu_0(\bfx)|^2 \; \rd \bfx.
\end{equation}

Finally, we have arrived at the bounds over time ($\mu \geq 0$):
\begin{equation}
\bfu,\; \upomega \in L^{\infty}(\Upo_T) \times  L^2(\rr) \cap C^2_0(\rr),\;\;\;
 \nabla p (\bfx) \in C^1(\rr).
\end{equation}
The most important result is written as
\begin{equation} \label{st2}
 \|\upomega(\bfx,t)\|_{L^{\infty}(\Upo_T)\times L^{\infty}(\rr)} < \infty.
\end{equation}
Admittedly, the present analyses on vorticity dynamics are much more precise and complete when compared with an earlier paper of Lam (2013). Note that the revised bounds are valid from $t>0$, including the Leray-Hopf regularity interval.

\subsection*{Uniqueness}

Let $(\bfv,q)$ be a second solution of the Navier-Stokes equations in addition to $(\bfu,p)$. Denote $\bfw=\bfv-\bfu$, and $\nabla{\cdot}\bfw=\nabla{\cdot}\bfv-\nabla{\cdot}\bfu=0$. The equation governing $\bfw$ reads
\begin{equation*}
\partial_t \bfw + (\bfv{\cdot}\nabla)\bfw + (\bfw{\cdot}\nabla)\bfu = \nabla (p-q)/\rho + \nu \Delta \bfw.
\end{equation*}
We multiply this equation by $\bfw$ and integrate the result over space. Then we estimate the magnitude of third term
\begin{equation}
\int_{\rr}|\bfw\cdot (\bfw{\cdot}\nabla)\bfu| \; \rd \bfx \leq \|\nabla \bfu\|_{L^{\infty}}\|\bfw\|^2_{L^2},
\end{equation}
by (\ref{d8}). Since the enstrophy $\|\nabla{\times}\bfw\|_{L^2}$ is never negative, the energy is determined by
\begin{equation}
\frac{\rd}{\rd t}\|\bfw \|^2_{L^2} \; \rd \bfx \leq 2 \|\nabla \bfu(\cdot,t) \|_{L^{\infty}} \|\bfw\|^2_{L^2}.
\end{equation}
The sum of the energy shows that the difference satisfies
\begin{equation} 
\|\bfv(\cdot,t)-\bfu(\cdot,t) \|^2_{L^2} \leq  C\; \|\bfv_0(\cdot)-\bfu_0(\cdot) \|^2_{L^2}.
\end{equation}
If the initial data of the flows coincide, solution $\bfu(\bfx,t)$ agrees with $\bfv(\bfx,t)$ everywhere. In conclusion, given the complete knowledge about the vorticity, the solutions of the Navier-Stokes equations are unique for any localised initial data of finite energy. A corollary proposition is that the Cauchy solutions of Euler's equations are well-posed for $t>0$.

It follows that the Cauchy problem of the incompressible vorticity equation possesses classical solutions for given smooth initial data having compact support. The solutions of the integral equation,
\begin{equation} \label{vts}
\begin{aligned}
\upomega(\bfx,t)=\int_{\rr} \! & H(\bfx{-}\bfy,t) \upomega_0(\bfy)\; \rd \bfy \\ 
&+ \int_0^t \! \int_{\rr} \! H(\bfx{-}\bfy,t{-}s) \;\big( (\upomega \cdot \nabla) \bfu  - (\bfu \cdot \nabla )\upomega\big) (\bfy,s) \; \rd \bfy \rd s,
\end{aligned}
\end{equation}
where
\begin{equation}
\bfu(\bfx) = -\frac{1}{4 \pi}\int_{\rr}\frac{\bfx-\bfy}{|\bfx-\bfy|^3}\times \upomega(\bfy) \; \rd \bfy,
\end{equation}
and
\begin{equation} \label{vt4}
\nabla\bfu(\bfx) = \frac{1}{4 \pi}\int_{\rr}\nabla\Big(\frac{1}{|\bfx-\bfy|}\Big)\nabla{\times}\upomega(\bfy) \; \rd \bfy,
\end{equation}
have been found to be a superposition of multiple convolutions of the initial data. 
The continuum dynamics instigates vortices or eddies of spatio-temporal scales of numerous ranges, see equation (8.11) of Lam (2013). The non-linearity in (\ref{vts}) is solely responsible for the production of the multitudinous structures, giving rise to turbulence. Critically, the field of the velocity-vorticity is self-consistent, and is kinematically regulated by the Helmholtz decomposition. The integral (\ref{vt4}) for the velocity gradient {\it does not contain a singular kernel}. In practice, it is almost unavoidable to involve heavy numerical computations in fluid simulation. 

\section{Conclusion}

The Cauchy problem of the primitive equations does not have any unique solution beyond the local regularity time. We have shown that the primitive setting for incompressible flow is over-specified in the sense that the dynamics can be split into a solenoidal diffusion and a reactive non-linear convection. The exception is the flows of zero pressure gradient. Our analysis clarifies the non-existence of strong {\it a priori} bounds, that are crucial in verifying the theoretical consistency. In essence, the Navier-Stokes equations are {\it incompletely} formulated. As asserted, the concept of unstable flow stems from the approximations due to the {\it ad hoc} linearisation. Above all, there exist no immediate instability counterparts in the exact vorticity dynamics. One consolation is that we understand why we cannot satisfactorily explain most of the well-established experimental observations about fluid flows, should we squarely cling on to the primitive perspectives.\footnote{The author recalled that, in an introductory undergraduate course on fluid mechanics, one of his professors said the question whether the velocity is induced by the pressure or vice versa is a chicken-and-egg argument.}

Fortunately, the global regularity of the vorticity equation for incompressible flows can be proved beyond doubt, thanks to the invariance of the vorticity field. 

\addcontentsline{toc}{section}{\noindent{References}}

\vspace{1cm}
\begin{acknowledgements}
\noindent 
28 March 2021

\noindent 
\texttt{f.lam11@yahoo.com}
\end{acknowledgements}
\end{document}